\documentclass[11pt, reqno]{amsart}
\usepackage{amssymb, amsmath, amsthm, mathrsfs}
\usepackage{hyperref}
\usepackage{enumitem}
\usepackage{float}
\newtheorem*{theo11A}{Theorem 1.1.A}
\newtheorem*{theo11B}{Theorem 1.1.B}

\newtheorem*{rem11A}{Remark 1.1.A}

\newtheorem*{theo12A}{Theorem 1.2.A}
\newtheorem*{theo12B}{Theorem 1.2.B}
\newtheorem*{theo12C}{Theorem 1.2.C}
\newtheorem*{theo12D}{Theorem 1.2.D}
\newtheorem*{theo12E}{Theorem 1.2.E}
\newtheorem*{theo12F}{Theorem 1.2.F}

\newtheorem*{cor12A}{Corollary 1.2.A}
\newtheorem*{cor12B}{Corollary 1.2.B}
\newtheorem*{exm12A}{Example 1.2.A}

\newtheorem{ques}{Question}[section]

\newtheorem{theo}{Theorem}[section]
\newtheorem{lem}{Lemma}[section]
\newtheorem{cor}{Corollary}[section]

\newtheorem{exm}{Example}[section]

\newtheorem{rem}{Remark}[section]

\newcommand{\ol}{\overline}

\newcommand{\be}{\begin{equation}}
\newcommand{\ee}{\end{equation}}
\newcommand{\beas}{\begin{eqnarray*}}
\newcommand{\eeas}{\end{eqnarray*}}
\newcommand{\bea}{\begin{eqnarray}}
\newcommand{\eea}{\end{eqnarray}}

\numberwithin{equation}{section}
\begin{document}

\title[O\MakeLowercase {n the existence of entire solutions to a system of nonlinear}......]{\LARGE O\Large\MakeLowercase {n the existence of entire solutions to a system of nonlinear} F\MakeLowercase {ermat}-\MakeLowercase {Type Partial Differential-Difference equations}}

\date{}
\author[J. F. X\MakeLowercase{u}, S. M\MakeLowercase {ajumder} \MakeLowercase{and} D. P\MakeLowercase {ramanik}]{J\MakeLowercase{unfeng} X\MakeLowercase{u}, S\MakeLowercase {ujoy} M\MakeLowercase {ajumder}$^*$ \MakeLowercase{and} D\MakeLowercase {ebabrata} P\MakeLowercase {ramanik}}
\address{Department of Mathematics, Wuyi University, Jiangmen 529020, Guangdong, People's Republic of China.}
\email{xujunf@gmail.com}
\address{Department of Mathematics, Raiganj University, Raiganj, West Bengal-733134, India.}
\email{sm05math@gmail.com, sjm@raiganjuniversity.ac.in}
\address{Department of Mathematics, Raiganj University, Raiganj, West Bengal-733134, India.}
\email{debumath07@gmail.com}

\renewcommand{\thefootnote}{}
\footnote{*\emph{Corresponding Author}: Sujoy Majumder.}

\renewcommand{\thefootnote}{\arabic{footnote}}
\setcounter{footnote}{0}

\begin{abstract} The aim of this study is to investigate the precise form of finite-order entire solutions to the following system of Fermat-type partial differential-difference equations:
\beas
\begin{cases}
\left(\frac{\partial f_1\left(z_1, z_2, \ldots, z_m \right)}{\partial z_1}\right)^{n_1} + f_2^{m_1} \left(z_1 + c_1, z_2 + c_2, \ldots, z_m + c_m \right) = 1,\\
\left(\frac{\partial f_2\left(z_1, z_2, \ldots, z_m \right)}{\partial z_1}\right)^{n_2} + f_1^{m_2} \left(z_1 + c_1, z_2 + c_2, \ldots, z_m + c_m \right) = 1
\end{cases}
\eeas

for various combinations of the positive integers $n_1$, $n_2$, $m_1$ and $m_2$. Our results extend the work of Xu et al. (Entire solutions for several systems of non-linear difference and partial differential-difference equations of Fermat-type, J. Math. Anal. Appl., 483(2), 2020), generalizing the setting $\mathbb{C}^2$ to $\mathbb{C}^m$. Several examples are provided to illustrate the applicability and sharpness of the obtained results.

\end{abstract}

\thanks{Typeset by \AmS -\LaTeX}
\maketitle

\tableofcontents

\section{{\bf Introduction and main results}}

\subsection{{\bf Fermat type functional equation in $\mathbb{C}$}} The following functional equation
\bea\label{ne1} f^n(z)+g^n(z)=1,\eea
where $n$ is a positive integer can be regarded as the Fermat diophantine equations $x^n+y^n=1$ over function fields.
In 1927, Montel \cite{M1} proved that the equation (\ref{ne1}) has no transcendental entire solutions for $n \geq 3$.
Gross \cite[Theorem 4]{FG2} proved that equation (\ref{ne1}) has no non-constant entire solutions when $n>2$. The equation (\ref{ne1}) has no non-constant meromorphic solutions when $n>3$ (see \cite[Theorem 3]{FG1}). For $n=2$, Gross \cite[Theorem 4]{FG2} found that the equation (\ref{ne1}) has entire solutions of the from $f(z)=\sin (h(z))$ and $g(z)=\cos (h(z))$, where $h(z)$ is an entire function.
In Theorem 1, Gross \cite{FG1} proved that all meromorphic solutions of the equation $f^2(z)+g^2(z)=1$ are the form
\[f(z)=\frac{1-\alpha^2(z)}{1+\alpha^2(z)}\;\;\text{and}\;\;g(z)=\frac{2 \alpha(z)}{1+\alpha^2(z)},\]
where $\alpha(z)$ is a non-constant meromorphic function. For $n=3$, Baker \cite[Theorem 1]{Baker} proved that the only non-constant meromorphic solutions of equation (\ref{ne1}) are the functions
\[f=\frac{1}{2}\left\{1+\frac{\wp^{(1)}(u)}{\sqrt{3}}\right\} / \wp(u)\;\;\text{and}\;\; g=\frac{\eta}{2}\left\{1-\frac{\wp^{(1)}(u)}{\sqrt{3}}\right\} / \wp(u)\]
for a non-constant entire function $u$ and a cubic root $\eta$ of unity, where
the Weierstrass $\wp$-function with periods $\omega_1$ and $\omega_2$ defined as
\[\wp(z;\omega_1,\omega_2)=\frac{1}{z^2}+\sum\limits_{\alpha,\beta:\alpha^2+\beta^2=1}\left\lbrace\frac{1}{(z+\alpha\omega_1+\beta\omega_2)^2}-\frac{1}{(\alpha\omega_1+\beta\omega_2)^2}\right\rbrace\]
and it satisfies $(\wp^{(1)})^2=4\wp^3-1$.

\smallskip
In 1970, Yang \cite{Y1} investigated the following Fermat-type equation
\bea\label{ne2} f^{m_1}(z)+g^{m_2}(z)=1\eea
and obtained that equation (\ref{ne2}) has no non-constant entire solutions when $\frac{1}{m_1}+\frac{1}{m_2}<1$ (see the proof of Theorem 1). Therefore it is clear that the equation (\ref{ne2}) has no non-constant entire solutions when $m_1>2$ and $m_2 > 2$. However, for the case when $m_1=m_2=2$ and $g(z)$ has a specific relationship with $f(z)$, many authors investigated the existence of solutions of the equation (\ref{ne2}). As a result, successively several research papers were published (see \cite{L.G, Liu1,LYL1, LCC1,LY1,Lu2,YL1}).

\smallskip
In 2012, Liu et al. \cite{LCC1} studied the existence of solutions for the differential-difference equation
\bea\label{ne3} (f^{(1)}(z))^{m_1} + f^{m_2}(z+c)=1\eea
and demonstrated that the equation (\ref{ne3}) has no finite order transcendental entire solutions when $m_1 \neq m_2$, where $m_1$ and $m_2$ are positive integers, and $c(\neq 0)$ is a constant. Additionally, Liu et al. \cite{LCC1} obtained the following result
\begin{theo11A}\cite[Theorem 1.3]{LCC1} The transcendental entire solutions with finite order of
\bea\label{sp1}(f^{(1)}(z))^2 + f^2(z+c) = 1\eea
must satisfy $f(z) = \sin (z \pm B \iota)$, where $B$ is a constant and $c = 2k\pi$ or $c = (2k+1)\pi$, with $k$ being an integer.
\end{theo11A}

\smallskip
Recall that the pair $(f(z), g(z))$ is referred to as a set of finite-order transcendental entire solutions for the system
\[\left\{\begin{array}{l}
f^{n_1}(z)+ g^{m_1}(z)=1 \\
f^{n_2}(z)+ g^{m_2}(z)=1
\end{array}\right.\]
if $f(z)$ and $g(z)$ are transcendental entire functions and $\rho=\max \{\rho(f), \rho(g)\}<\infty$.

\smallskip
In 2016, Gao \cite{L.G} addressed the problem of determining the form of solutions for the following system of Fermat-type differential-difference equations
\bea\label{sp}
\left\{\begin{array}{l}
{\left(f_1^{(1)}(z)\right)^2+{f_2}^2(z+c)^2=1} \\
{\left(f_2^{(1)}(z)\right)^2+{f_1}^2(z+c)=1}
\end{array}\right.\eea
and derived the following result:

\begin{theo11B}\cite[Theorem 1.1]{L.G}
Suppose that $(f_1(z), f_2(z))$ is a pair of finite order transcendental entire solutions for the system of differential-difference equations (\ref{sp}). Then $(f_1(z), f_2(z))$ satisfies one of the following
\begin{enumerate}
\item[(i)] $(f_1(z), f_2(z))=\left(\sin (z-b\iota), \sin (z-b_1 \iota)\right)$,
\item[(ii)] $(f_1(z), f_2(z))=\left(\sin (z+b\iota), \sin (z+b_1 \iota)\right)$,
\end{enumerate}
where $b, b_1$ are constants, and $c=k \pi$, $k$ is an integer.
\end{theo11B}

\begin{rem11A} When $f_1(z)\equiv f_2(z)\equiv f(z)$, it is easy to verify that the system of equations (\ref{sp}) reduces to the equation (\ref{sp1}). Now, from the conclusion of Theorem 1.1.B, it is easy to get: $f(z)=\sin (z\pm b\iota)$, where $b$ is a constant, and $c=k \pi$, $k$ is an integer. This shows that the conclusions of Theorem 1.1.B contain the conclusion of Theorem 1.1.A.
\end{rem11A}

\subsection{{\bf Fermat type functional equation in $\mathbb{C}^m$}}

The following theorem describes the entire and meromorphic solutions in $\mathbb{C}^m$ of the Fermat-type Eq.
\bea\label{FT1} f^n+g^n=1,\;n>1.\eea

\begin{theo12A} \cite[Theorem 1.3]{GS2} For $h:\mathbb{C}^m \to \mathbb{C}$ entire, the solutions of the Eq. (\ref{FT1}) are characterized as follows:
\begin{enumerate}
\item[(a)] for $n=2$, the entire solutions are $f=\cos(h)$ and $g=\sin(h)$;
\item[(b)] for $n>2$, there are no non-constant entire solutions;
\item[(c)] for $n=2$, the meromorphic solutions are of the form
$f=\frac{1 - \beta^2}{1 + \beta^2}$ and $g=\frac{2\beta}{1 + \beta^2}$,
with $\beta$ being meromorphic on $\mathbb{C}^m$;
\item[(d)] for $n=3$, the meromorphic solutions are of the form
\[f = \frac{1}{2\wp(h)}\left( 1 + \frac{\wp^{(1)}(h)}{\sqrt{3}}\right)\;\;\text{and}\;\;g = \frac{1}{2\wp(h)} \left(1-\frac{\wp^{(1)}(h)}{\sqrt{3}}\right);\]
\item[(e)] for $n >3$, there are no non-constant meromorphic solutions.
\end{enumerate}
\end{theo12A}

\smallskip
We define $\mathbb{Z}_+=\mathbb{Z}[0,+\infty)=\{n\in \mathbb{Z}: 0\leq n<+\infty\}$ and $\mathbb{Z}^+=\mathbb{Z}(0,+\infty)=\{n\in \mathbb{Z}: 0<n<+\infty\}$. We define
\[\displaystyle \partial_{z_i}(f(z))=\frac{\partial f(z)}{\partial z_i},\ldots,\partial^{l_i}_{z_i}(f(z))=\frac{\partial^{l_i} f(z)}{\partial z_i^{l_i}},\;\;\partial^{I}(f(z))=\frac{\partial^{|I|}f(z)}{\partial z_1^{i_1}\cdots \partial z_m^{i_m}}\]
where $l_i\in \mathbb{Z}_+$, $i=1,2,\ldots,m$ and $I=(i_1,i_2,\ldots,i_m)\in\mathbb{Z}^m_+$ such that $|I|=\sum_{j=1}^m i_j$.

\smallskip
In 2008, Li \cite[Theorem 1]{L1} proved that meromorphic solutions of the Fermat-type functional equation $f^2+g^2=1$ in $\mathbb{C}^2$ must be constant if and only if $\partial_{z_2} (f(z_1,z_2))$ and $\partial_{z_1} (g(z_1,z_2))$ share $0$ CM.  Moreover if $f(z_1,z_2)=\partial_{z_1}(u(z_1,z_2))$ and $g(z_1,z_2)=\partial_{z_2}(u(z_1,z_2))$, then any entire solutions of the partial differential equations:
\[\left(\partial_{z_1}(u(z_1,z_2))\right)^2+\left(\partial_{z_2}(u(z_1,z_2))\right)^2=1\]
in $\mathbb{C}^2$ are necessarily linear (see \cite{DK1}). From the fact that the surface $z_1^n+z_2^n= 1$ in $\mathbb{C}^2$ is a Kobayashi hyperbolic manifold and thus there are no non-constant holomorphic curves $(f_1, f_2)$ from the complex
plane $\mathbb{C}$ (and thus $\mathbb{C}^m$) to the surface (see \cite[page 360]{BS1}), one can easily prove that any
entire solutions of the partial differential equations
\[\left(\partial_{z_1}(u(z_1,z_2))\right)^n+\left(\partial_{z_2}(u(z_1,z_2))\right)^n=1,\]
where $n\geq 3$ must also be linear. The solutions of Fermat-type partial differential equation 
was originally investigated by Li \cite{L2} and Saleeby \cite{GS1}.

\smallskip
Recently, many authors  (see \cite{AA1}-\cite{AM1}, \cite{MSP}, \cite{XC1}, \cite{XLL1}, \cite {XW1}, \cite{XH}) have paid their attention to explore the existence of solutions of Fermat-type difference and partial differential-difference equations using the difference Nevanlinna theory for meromorphic functions (see \cite{TBC1, CK1, CX1, K1}) of higher dimension.

\smallskip
In 2018, Xu and Cao \cite{XC1} investigated the existence of the finite order transcendental entire solutions for a specific type of  Fermat-type partial differential-difference equations and obtained the following results.

\begin{theo12B}\cite[Theorem 1.1]{XC1} Let $c=(c_1, c_2)\in \mathbb{C}^2$. Then the Fermat-type partial differential-difference equation
\bea\label{f01}\left(\partial_{z_1}(f(z_1, z_2))\right)^{m_1}+f^{m_2}(z_1+c_1, z_2+c_2)=1\eea
doesn't have any transcendental entire solution with finite order, where $m_1$ and $m_2$ are two distinct positive integers.
\end{theo12B}

According to Theorem 1.2.B, the finite order transcendental entire solutions of (\ref{f01}) exist only possibly for the case $m_1=m_2$.
For finite order transcendental entire solutions of the following Fermat-type partial differential-difference equation
\bea\label{f1}\left(\partial_{z_1}(f(z_1, z_2))\right)^2+f^2(z_1+c_1, z_2+c_2)=1\eea
Xu and Cao \cite{XC2} obtained the following result.

\begin{theo12C}\cite[Theorem 1.1]{XC2} Let $c=(c_1, c_2)$ be a constant in $\mathbb{C}^2$. Then any finite order transcendental entire solutions of the Fermat-type partial difference-differential equation (\ref{f1})
has the form of $f(z_1, z_2)=\sin(A z_1+Bz_2+H(z_2))$,
where $A$ and $B$ are constants on $\mathbb{C}$ satisfying $A^2=1$ and $A e^{i (Ac_1+Bc_2)}=1$ and $H(z_2)$ is a polynomial in one variable $z_2$ such that $H(z_2)=H(z_2+c_2)$. In the special case whenever $c_2=0$, we have
$f(z_1, z_2)=\sin (Az_1+Bz_2+\text{Constant})$.
\end{theo12C}

If $c=(0, 0)$, then Theorem 1.2.C implies the following corollary.
\begin{cor12A}\cite[Corollary 1.4]{XC2} Any transcendental entire solution with finite order of the partial differential equation of the Fermat-type
\beas\left(\partial_{z_1}(f(z_1, z_2))\right)^2+f^2(z_1, z_2)=1\eeas
has the form of $f(z_1,z_2)=\sin (z_1+g(z_2))$, where $g(z_2)$ is a polynomial in one variable $z_2$.
\end{cor12A}

Also in the same paper, Xu and Cao \cite{XC2} exhibited the following example to explain the conclusion of Theorem 1.2.C.

\begin{exm12A}\cite[Example 1.2.]{XC2} Let $A=1$, $B=2$, and let two constants $c_1$ and $c_2$ satisfy $e^{\iota c_1}=1$ and $c_2=0$. Then $A e^{\iota (Ac_1+Bc_2)}=1$. The entire function $f(z)=\sin \left(z_1+2z_2+z_2^3+1\right)$ satisfies the following Fermat-type partial differential-difference equation
\beas \left(\partial_{z_1}(f(z_1, z_2))\right)^2+f^2(z_1+c_1, z_2+c_2)=1\;\text{in}\;\mathbb{C}^2,\; \text{where}\;c=(c_1, c_2).\eeas
\end{exm12A}

Example 1.2.A shows that the function $H(z_2)$ in the conclusion of Theorem 1.2.C may be a non-constant polynomial whenever $c_2=0$.

\medskip
In 2020, Xu and Wang \cite {XW1} again considered Eq. (\ref{f01}) to discuss the existence of finite order transcendental entire solutions for the case when $m_1\neq m_2$. Actually Xu and Wang \cite[Remark 1.7]{XW1} proved that: If one of the following conditions holds:
\begin{enumerate}
\item[(i)] $m_2>m_1$,
\item[(ii)] $m_1>m_2\geq 2$,
\end{enumerate}
then Eq. (\ref{f01}) does not have any finite order transcendental entire solutions.

\smallskip
When $m_1=2$ and $m_2=1$,  Xu and Wang \cite{XW1} obtained the following result regarding the finite order transcendental entire solutions of Eq. (\ref{f01}).

\begin{theo12D}\cite[see pp.8]{XW1} Let $c=(c_1, c_2)\in\mathbb{C}^2$. If $f(z_1,z_2)$ is a finite order transcendental entire solution of the following Eq.
\bea\label{f03} \left(\partial_{z_1}(f(z_1, z_2))\right)^{2}+f(z_1+c_1, z_2+c_2)=1,\eea
then $f(z_1,z_2)$ takes the from
\beas &&f(z_1,z_2)\\&=&1-\frac{1}{4}c_1^2-\frac{1}{4}z_1^2+\frac{c_1}{2c_2}z_1z_2-\frac{c_1^2}{2c_2}(z_2-c_2)+(z_1-c_1)G_1(z_2)-\left\lbrack \frac{c_1}{2c_2}(z_2-c_2)+G_1(z_2)\right\rbrack^2,\eeas
where $G_1(z_2)$ is a finite order transcendental entire period function with period $c_2$.
\end{theo12D}

\medskip
Recently Xu et al. \cite{XLL1} studied the existence and the forms of the finite order transcendental entire solutions to the following system of Fermat-type partial differential-difference equations:
\bea\label{xll}
\begin{cases}
\left(\partial_{z_1}(f_1(z_1, z_2))\right)^{n_1} + f_2^{m_1}(z_1 + c_1, z_2 + c_2) = 1,\\
\left(\partial_{z_1}(f_2(z_1, z_2))\right)^{n_2} + f_1^{m_2}(z_1 + c_1, z_2 + c_2) = 1,
\end{cases}
\eea
where $m_i$ and $n_i$ are positive integers for $i=1,2$.

\medskip
Firstly, in the following result, Xu et al. \cite{XLL1} proved that the system of Fermat-type partial differential-difference equations (\ref{xll}) does not have any pair of finite order transcendental entire solutions under some restrictions on $n_i$ and $m_i$ for $i=1,2$.

\begin{theo12E} \cite[Theorem 1.1]{XLL1} Let $c=(c_1, c_2)\in\mathbb{C}^2$ and let $m_i$ and $n_i$ ($i=1,2$) be positive integers. If the system of Fermat-type partial differential-difference equations (\ref{xll}) satisfies one of the following conditions:
\begin{enumerate}
\item[(i)] $m_1 m_2>n_1 n_2$;
\item[(ii)] $m_j>\frac{n_j}{n_j-1}$ for $n_j \geq 2$, $j=1,2$,
\end{enumerate}
then the system of equations (\ref{xll}) does not have any pair of finite order transcendental entire solutions.
\end{theo12E}

\medskip
Secondly, in the following result, Xu et al. \cite{XLL1} derived the solutions of the system of Fermat-type differential-difference equations (\ref{xll}) for the case when $n_1=n_2=m_1=m_2=2$.

\begin{theo12F} \cite[Theorem 1.3]{XLL1} Let $c=(c_1, c_2)\in \mathbb{C}^2$. Then any pair of finite order transcendental entire solutions for the following system of Fermat-type partial differential-difference equations
\bea\label{xll1}
\begin{cases}
\left(\partial_{z_1}(f_1(z_1, z_2))\right)^2+f_2^2(z_1 + c_1, z_2 + c_2) = 1,\\
\left(\partial_{z_1}(f_2(z_1, z_2))\right)^2+f_1^2(z_1 + c_1, z_2 + c_2) = 1,
\end{cases}
\eea
have the following forms
\[\left(f_1(z), f_2(z)\right)=\left(\frac{e^{L(z)+B_1}+e^{-\left(L(z)+B_1\right)}}{2}, \frac{A_{21} e^{L(z)+B_1}+A_{22} e^{-\left(L(z)+B_1\right)}}{2}\right),\]
where $L(z)=a_1 z_1+a_2 z_2$, $B_1\in\mathbb{C}$ and $a_1, c, A_{21}, A_{22}$ satisfy one of the following cases
\begin{enumerate}
\item[(i)] $A_{21}=-\iota, A_{22}=\iota$, and $a_1=\iota, L(c)=\left(2 k+\frac{1}{2}\right) \pi \iota$, or $a_1=-\iota, L(c)=\left(2 k-\frac{1}{2}\right) \pi \iota$;
\item[(ii)] $A_{21}=\iota, A_{22}=-\iota$, and $a_1=\iota, L(c)=\left(2 k-\frac{1}{2}\right) \pi \iota$, or $a_1=-\iota, L(c)=\left(2 k+\frac{1}{2}\right) \pi \iota$;
\item[(iii)] $A_{21}=1, A_{22}=1$, and $a_1=\iota, L(c)=2 k \pi \iota$, or $a_1=-\iota, L(c)=(2 k+1) \pi \iota$;
\item[(iv)] $A_{21}=-1, A_{22}=-1$, and $a_1=\iota, L(c)=(2 k+1) \pi \iota$, or $a_1=-\iota, L(c)=2 k \pi \iota$.
\end{enumerate}
\end{theo12F}

If $c=(0,0)$, then Theorem 1.2.F implies the following corollary:
\begin{cor12B} \cite[Corollary 1.1]{XLL1} Any pair of finite order transcendental entire solutions for the following system of Fermat-type partial differential equations
\beas
\begin{cases}
\left(\partial_{z_1}(f_1(z_1, z_2))\right)^2+f_2^2(z_1, z_2) = 1,\\
\left(\partial_{z_1}(f_2(z_1, z_2))\right)^2+f_1^2(z_1, z_2) = 1,
\end{cases}
\eeas
have the following forms
\[\left(f_1(z), f_2(z)\right)=\left(\frac{e^{L(z)+B_1}+e^{-\left(L(z)+B_1\right)}}{2}, \eta\frac{e^{L(z)+B_1}+e^{-\left(L(z)+B_1\right)}}{2}\right),\]
where $L(z)=a_1 z_1+a_2 z_2$, $B_1\in\mathbb{C}$ and $a_1, \eta$ satisfy $(i)$ $a_1=\iota$ and $\eta=1$; $(ii)$ $a_1=-\iota$ and $\eta=-1$.
\end{cor12B}

\smallskip
In the study of existence and the forms of the finite order transcendental entire solutions to the system of the Fermat-type partial differential-difference equations (\ref{xll}), Theorems 1.2.E and 1.2.F unquestionably play a significant role. But, if we take $f_1(z)\equiv f_2(z)\equiv f(z)$, then the system of equations (\ref{xll1}) reduces to the equation (\ref{f1}).
Therefore from the conclusions of Theorem 1.2.F, it is easy to get: $f(z)=\sin (A z_1+B z_2+B_1)$, where $A$, $B$ and $B_1$ are constants in $\mathbb{C}$ satisfying $A^2=1$ and $A e^{i (Ac_1+Bc_2)}=1$. One may therefore conclude that the result of Theorem 1.2.F for the case where $f_1(z)\equiv f_2(z)$ is not complete by looking at Example 1.2.A.

Next, by considering the following example something more can be said about the precise form of the finite order transcendental entire solutions for the system of equations (\ref{xll1}) than that are obtained in Theorem 1.2.F.

\begin{exm}\label{exm11} Let
\[(f_1(z),f_2(z))=(\sin(Az_1+A_{12}z_2+Q_1(z_2)), \sin(Bz_1+B_{12}z_2+Q_2(z_2)),\]
where $A=A_{12}=B=B_{12}=1$, $Q_1(z_2)=z_2^2$, $Q_2(z_2)=z_2^2+\pi$ and $(c_1,c_2)=(2\pi,0)$. Note that $Be^{2\iota(Ac_1+A_{12}c_2)}=A$ and $Ae^{2\iota(Bc_1+B_{12}c_2)}=B$. It is easy to verify that $(f_1(z), f_2(z))$ satisfy the following system of equations
\beas
\begin{cases}
\left(\partial_{z_1}(f_1(z_1, z_2))\right)^2+f_2^2(z_1+2\pi, z_2) = 1,\\
\left(\partial_{z_1}(f_2(z_1, z_2))\right)^2+f_1^2(z_1+2\pi, z_2) = 1.
\end{cases}
\eeas
\end{exm}

The conclusion of Theorem 1.2.F shows that the order of the transcendental entire solutions $(f_1(z),f_2(z))$ is $1$, i.e, $\rho(f_i)=1$ for $i=1,2$. However, Example \ref{exm11} demonstrates that the system of equations (\ref{xll1}) may have transcendental entire solutions $(f_1(z),f_2(z))$ with an order greater than $1$.

\smallskip
Finding the exact form of the finite order transcendental entire solutions $(f_1(z),f_2(z))$ of the system of equations (\ref{xll1}) would therefore be intriguing. We also want to generalize Theorem 1.2.F in this study by taking into account the following system of  Fermat-type partial differential-difference equations
\bea\label{pds1}
\begin{cases}
\left(\partial_{z_1}(f_1(z_1, z_2, \ldots, z_m))\right)^{n_1} + f_2^{m_1}(z_1 + c_1, z_2 + c_2, \ldots, z_m + c_m) = 1, \\

\medskip
\left(\partial_{z_1}(f_2(z_1, z_2, \ldots, z_m))\right)^{n_2} + f_1^{m_2}(z_1 + c_1, z_2 + c_2, \ldots, z_m + c_m) = 1
\end{cases}
\eea
in $\mathbb{C}^m$, where $n_1, n_2, m_1, m_2$ are positive integers, and $c_1, c_2, \ldots, c_m$ are constants in $\mathbb{C}$.
For $n_1=n_2=m_1=m_2=1$, it is easy to verify that $(f_1(z),f_2(z))=(e^{z_1+z_2+\ldots+z_m}+1,e^{z_1+z_2+\ldots+z_m}+1)$ is a solution of the Eq. (\ref{pds1}), where $e^{c_1+c_2+\ldots+c_m}=-1$. Therefore we consider the Eq. (\ref{pds1}) for the
existence of entire solutions for the case when $n_i+m_i>2$, where $i=1,2$.

\medskip
In this paper, our next goal is to investigate all possible finite-order transcendental entire solutions $(f_1(z),f_2(z))$ for the system of equations (\ref{pds1}) for all different combinations of $n_1$, $n_2$, $m_1$ and $m_2$. Specifically, we focus on understanding the conditions under which such solutions exist and determining their explicit forms.  
Furthermore, we extend these findings to higher-dimensional settings $\mathbb{C}^m$ and explore how the solutions behave when multiple variables are involved. Now we state our main results.

\begin{theo}\label{t2.1} Let $c\in\mathbb{C}^m$. Then any finite order transcendental entire solutions for the system of equations (\ref{pds1}), where $m_1,m_2,n_1,n_2\in\mathbb{N}$ such that $m_i+n_i>2$, $i=1,2$ are characterized as follows:
\begin{enumerate}
\item[(1)] if $n_1=m_1$, then $n_1=n_2=m_1=m_2=2$ and
\[(f_1(z),f_2(z))=\sin (L_1(z)+Q_1(z_2,\ldots, z_m)),\sin (L_2(z)+Q_2(z_2,\ldots, z_m))),\]
where $L_1(z)=Az_1+A_{11}z_2\ldots+A_{1m}z_m$, $L_2(z)=Bz_1+B_{11}z_2\ldots+B_{1m}z_m$, $Q_i(z_2,\ldots,z_m)$ is a polynomial such that $Q_i(z_2+2c_2,\ldots,z_m+2c_m)=Q_i(z_2,\ldots,z_m)$ for $i=1,2$, $A_{1i},B_{1i}\in\mathbb{C}$ for $i=2,\ldots,m$ such that $A^2=1$, $B^2=1$ and one of the following hold:
\begin{enumerate}
\item[(i)] $Be^{2\iota (Ac_1+A_{12}c_2+\ldots+A_{1m}c_m)}=A$ and $Ae^{2\iota (Bc_1+B_{12}c_2+\ldots+B_{1m}c_m)}=B$;

\smallskip
\item[(ii)] $ABe^{2\iota (Ac_1+A_{12}c_2+\ldots+A_{1m}c_m)}=1$ and $ABe^{2\iota (Bc_1+B_{12}c_2+\ldots+B_{1m}c_m)}=1$;
\end{enumerate}
\item[(2)] if one of the following conditions hold:
\begin{enumerate}
\item[(2)(i)] $m_1m_2>n_1n_2$,

\smallskip
\item[(2)(ii)] $n_i=m_i$ and $n_j>m_j$ for $i,j\in\{1,2\}$ such that $i\neq j$,

\smallskip
\item[(2)(iii)] $n_i>m_i$ and $n_1n_2-n_i>2$ for $i\in\{1,2\}$,

\smallskip
\item[(2)(iv)] $n_i>m_i$ and $m_j>n_j$ such that $n_i\geq 3$ and $m_j>\frac{n_i}{n_i-2}$, where $i\neq j$ and $i,j\in\{1,2\}$,
\end{enumerate}
then the system of equations (\ref{pds1}) doesn't have any finite order transcendental entire solutions;
\item[(3)] if $n_i>m_i$, then $m_i=1$ and $n_i=2$  for $i=1,2$ and
\beas f_i(z)=1+\frac{K_i}{4}z_1^2 + z_1g_i(z_2,z_3,\ldots,z_m)-K_i^2g_i^2(z_2,z_3,\ldots,z_m)\eeas
where $K_i^3=-1$ and $g_i(z_2,z_3,\ldots, z_m)$ is a finite order transcendental entire function such that
$g_i(z_2+2c_2, \ldots, z_m+2c_m)-g_i(z_2, \ldots, z_m)=K_ic_1$ for $i=1,2$.
\end{enumerate}
\end{theo}

We give two examples to explain the conclusion $(1)$ of Theorem \ref{t2.1}.
\begin{exm} Suppose that
\[(f_1(z),f_2(z))=(\sin(Az_1+A_{12}z_2+A_{13}z_3+Q_1(z_2,z_3)), \sin(Bz_1+B_{12}z_2+B_{13}z_3+Q_2(z_2,z_3)),\]
where $A=A_{12}=A_{13}=B=B_{12}=B_{13}=1$, $Q_1(z_2,z_3)=(z_2-z_3)^2$, $Q_2(z_2,z_3)=(z_2-z_3)^2+\pi$ and $(c_1,c_2,c_3)=(0, \frac{\pi}{2}, \frac{\pi}{2})$. Note that $Be^{2\iota(Ac_1+A_{12}c_2+A_{13}c_3)}=A$ and $Ae^{2\iota(Bc_1+B_{12}c_2+B_{13}c_3)}=B$. It is easy to verify that $(f_1(z), f_2(z))$ satisfy the following system of equations
\bea\label{pds2}
\begin{cases}
\left(\partial_{z_1}(f_1(z_1,z_2,z_3))\right)^2+f_2^2(z_1+c_1,z_2+c_2,z_3+c_3)=1, \\
\left(\partial_{z_1}(f_2(z_1,z_2,z_3))\right)^2+f_1^2(z_1+c_1,z_2+c_2,z_3+c_3)=1.
\end{cases}
\eea
\end{exm}

\begin{exm} Suppose that
\[(f_1(z),f_2(z))=(\sin(Az_1+A_{12}z_2+A_{13}z_3+Q_1(z_2,z_3)), \sin(Bz_1+B_{12}z_2+B_{13}z_3+Q_2(z_2,z_3))),\]
where $A=A_{12}=1$, $A_{13}=-1$, $B=B_{12}=1$, $B_{13}=-1$, $Q_1(z_2,z_3)=Q_2(z_2,z_3)=(z_2-z_3)^3$ and $(c_1,c_2,c_3)=(0, \pi, \pi)$. Note that $ABe^{\iota(Ac_1+A_{12}c_2+A_{13}c_3)}=1$ and $ABe^{2\iota(Bc_1+B_{12}c_2+B_{13}c_3)}=1$. It is easy to verify that $(f_1(z), f_2(z))$ satisfy the system of equations (\ref{pds2}).
\end{exm}

Following example is provided to illustrate our conclusion $(3)$ of Theorem \ref{t2.1}.

\begin{exm} Suppose that
\[f_i(z_1, z_2, z_3)=1-\frac{1}{4} z_1^2 +g_i(z_2,z_3)z_1 -g_i^2(z_2,z_3),\]
where $g_i(z_2,z_3)=\sin (z_2 + z_3)$ and $c=(0, \pi, \pi)$. We take $K_i=-1$. Note that $g_i(z_2,z_3)$ is a transcendental entire function of finite order such that $g_i(z_2+2c_2,z_3+2c_3)-g_i(z_2,z_3)=0=K_ic_1$. It is easy to verify that $(f_1(z),f_2(z))$ satisfy the following system of equations
\bea\label{pds2a}
\begin{cases}
\left(\partial_{z_1}(f_1(z_1,z_2,z_3))\right)^2+f_2(z_1+c_1,z_2+c_2,z_3+c_3)=1, \\
\left(\partial_{z_1}(f_2(z_1,z_2,z_3))\right)^2+f_1(z_1+c_1,z_2+c_2,z_3+c_3)=1.
\end{cases}
\eea
\end{exm}

When $f_1(z)\equiv f_2(z)\equiv f(z)$, $n_2=n_1$ and $m_2=m_1$, it is easy to see that system (\ref{pds1}) can be reduced as the following equation
\bea\label{sp4} \left(\partial_{z_1}(f(z_1, z_2, \ldots, z_m))\right)^{n_1} + f^{m_1}(z_1 + c_1, z_2 + c_2, \ldots, z_m + c_m) = 1\eea
and so by Theorem \ref{t2.1}, we get the following corollary:

\begin{cor}\label{c2.1} Let $c\in\mathbb{C}^m$. Then any finite order transcendental finite order entire solutions of the equation (\ref{sp4}), where $m_1,n_1\in\mathbb{N}$ such that $m_1+n_1>2$ are characterized as follows:
\begin{enumerate}
\item[(i)] if $n_1=m_1$, then $n_1=m_1=2$ and
\[f(z)=\sin (Az_1+\ldots+A_mz_m+P(z_2,\ldots, z_m)),\]
where $A$, $A_2,\ldots$, $A_m\in\mathbb{C}$ such that $A^2 = 1$ and $Ae^{\iota(Ac_1+\ldots+A_mc_m)}=1$ and $P(z_2,\ldots, z_m)$ is a polynomial such that $P(z_2+c_2,\ldots,z_m+c_m)\equiv P(z_2,\ldots,z_m)$. In the special case
whenever $m=2$ and $c_2\neq 0$, we get $f(z_1, z_2)=\sin(A_1z_1+A_2z_2+\text{Constant})$,

\smallskip
\item[(ii)] if $n_1\neq m_1$, then $n_1=2$, $m_1=1$ and
\[f(z)=1+\frac{K_i}{4} z_1^2 +z_1g(z_2, \ldots, z_m)-K_i^2 g^2(z_2, \ldots, z_m),\]
where $g(z_2,\ldots, z_m)$ is a transcendental entire function of finite order such that $g(z_2 + c_2, \ldots, z_m + c_m) - g(z_2, \ldots, z_m) = \frac{1}{2}c_1$.
\end{enumerate}
\end{cor}

\begin{rem} It is evident that Corollary \ref{c2.1} extends Theorem 1.2.C from $\mathbb{C}^2$ to $\mathbb{C}^m$.
\end{rem}

If  $c=(0,0,\ldots, 0)$, then Theorem \ref{t2.1} implies the following corollary.

\begin{cor}\label{c2.2} The transcendental entire solutions for the following system of Fermat-type partial differential equations
\bea\label{pds1a}
\begin{cases}
\left(\partial_{z_1}(f_1(z_1, z_2, \ldots, z_m))\right)^{n_1} + f_2^{m_1}(z_1, z_2, \ldots, z_m) = 1, \\
\left(\partial_{z_1}(f_2(z_1, z_2, \ldots, z_m))\right)^{n_2} + f_1^{m_2}(z_1, z_2, \ldots, z_m) = 1,
\end{cases}
\eea
where $m_1,m_2,n_1,n_2\in\mathbb{N}$ such that $m_i+n_i>2$, $i=1,2$ are characterized as follows:
\begin{enumerate}
\item[(1)] if $n_1=m_1$, then $n_1=n_2=m_1=m_2=2$ and
\[(f_1(z),f_2(z))=(\sin (L_1(z)+Q_1(z_2,\ldots, z_m)),\sin (L_2(z)+Q_2(z_2,\ldots, z_m))),\]
where $L_1(z)=Az_1+A_{11}z_2\ldots+A_{1m}z_m$, $L_2(z)=Bz_1+B_{11}z_2\ldots+B_{1m}z_m$, $Q_i(z_2,\ldots,z_m)$ is a polynomial for $i=1,2$, $A_{1i},B_{1i}\in\mathbb{C}$ for $i=2,\ldots,m$ such that either $A^2=B^2=1$ or $AB=1$;

\smallskip
\item[(2)] if one of the following conditions hold:
\begin{enumerate}
\item[(2)(i)] $m_1m_2>n_1n_2$,

\smallskip
\item[(2)(ii)] $n_i=m_i$ and $n_j> m_j$ for $i,j\in\{1,2\}$ such that $i\neq j$,

\smallskip
\item[(2)(iii)] $n_i>m_i$ and $n_1n_2-n_i>2$ for $i\in\{1,2\}$,

\smallskip
\item[(2)(iv)] $n_i>m_i$ and $m_j>n_j$ such that $n_i\geq 3$ and $m_j>\frac{n_i}{n_i-2}$, where $i\neq j$ and $i,j\in\{1,2\}$,
\end{enumerate}
then the system of equations (\ref{pds1a}) doesn't have any finite order transcendental entire solutions;
\item[(3)] if $n_i>m_i$, then $m_i=1$ and $n_i=2$  for $i=1,2$ and
\beas f_i(z)=1+\frac{K_i}{4}z_1^2 + z_1g_i(z_2,z_3,\ldots,z_m)-K_i^2g_i^2(z_2,z_3,\ldots,z_m)\eeas
where $K_i^3=-1$ and $g_i(z_2,z_3,\ldots, z_m)$ is a transcendental entire functions for $i=1,2$.
\end{enumerate}
\end{cor}

We now make the following observations on Theorem \ref{t2.1}.

\begin{rem}
We have explored as many scenarios as possible for all values of $n_i$ and $m_i$, where $i=1,2$. Theorems 1.2.E and 1.2.F did not address the following cases, which we have covered in Theorem \ref{t2.1} as follows:
\begin{enumerate}
\item[(i)] For $(n_1,m_1,n_2,m_2)=(2,2,2,1)$, $(n_1,m_1,n_2,m_2)=(2,1,2,2)$ and $(n_1,m_1,n_2,m_2)=(3,1,2,1)$,
we observed from part-$(2)$ of Theorem \ref{t2.1} that the system of equations (\ref{pds1}) does not admit any finite order transcendental entire solutions. Moreover, for $n_i>m_i$ and $m_j>n_j$ such that $n_i\geq 3$ and $m_j>\frac{n_i}{n_i-2}$, where $i,j\in\{1,2\}$ such that $i\neq j$, we demonstrated that the system of equations (\ref{pds1}) does not admit any finite order transcendental entire solutions.

\medskip
\item[(ii)] For $n_i>m_i\;(i=1,2)$, we established that the system of Fermat-type partial differential-difference equations (\ref{pds1}) have finite order transcendental entire solutions of the form
\[f_i(z) = 1 + \frac{K_i}{4}z_1^2 + z_1g_i(z_2, z_3, \ldots, z_m)-K_i^2 g_i^2(z_2, z_3, \ldots, z_m),\]
where \( K_i^3 = -1 \) and \( g_i(z_2, z_3, \ldots, z_m) \) is a finite order transcendental entire function satisfying
$g_i(z_2 + 2c_2, \ldots, z_m + 2c_m) - g_i(z_2, \ldots, z_m) = K_ic_1$ for $i=1,2$.
\end{enumerate}

This demonstrates unequivocally that by offering a thorough analysis of every additional scenario that could arise, our findings significantly outperform those of Xu et al. \cite{XLL1}.
\end{rem}

Now based on Theorem \ref{t2.1}, we ask the following question:

\begin{ques}
What can be inferred about the existence of finite order transcendental entire solutions for the system of Fermat-type partial differential equations (\ref{pds1}) without the hypothesis ``$n_i>m_i$ and $m_j>n_j$ such that $n_i\geq 3$ and $m_j>\frac{n_i}{n_i-2}$, $i,j\in\{1,2\}$ such that $i\neq j$''?
\end{ques}

\subsection{\bf Basic Notations in several complex variables}
As on $\mathbb{C}^m$, the exterior derivative $d$ splits $d= \partial+ \bar{\partial}$ and twists to $d^c= \frac{i}{4\pi}\left(\bar{\partial}- \partial\right)$ (see \cite{HLY1,WS1}). Clearly $dd^{c}= \frac{i}{2\pi}\partial\bar{\partial}$.
An exhaustion $\tau_m$ of $\mathbb{C}^m$ is defined by $\tau_m(z)=||z||^2$. The standard Kaehler metric on $\mathbb{C}^m$ is given by $\upsilon_m=dd^c\tau_m>0$. On $\mathbb{C}^m\backslash \{0\}$, we define $\omega_m=dd^c\log \tau_m\geq 0$ and $\sigma_m=d^c\log \tau_m \wedge \omega_m^{m-1}$. For any $S\subseteq \mathbb{C}^m$, let $S[r]$, $S(r)$ and $S\langle r\rangle$ be the intersection of $S$ with
respectively the closed ball, the open ball, the sphere of radius $r>0$ centered at $0\in \mathbb{C}^m$.

\smallskip
The zero multiplicity $\mu^0_f(a)$ of an entire function $f$ in $\mathbb{C}^m$ at a point $a\in \mathbb{C}^m$ is defined
to be the order of vanishing of $f$ at $a$. If $a=(a_1,\ldots,a_m)$, then $\partial^I(f(a))=0$, where $|I|\leq \mu^0_f(a)-1$.
In other words, we can write $f(z)=\sum_{i=0}^{\infty}P_i(z-a)$, where the term $P_i(z-a)$ is either identically zero or a homogeneous polynomial of degree $i$. Certainly $\mu^0_f(a)=\min\{i:P_i(z-a)\not\equiv 0\}$.

\smallskip
Let $f$ be a meromorphic function in $\mathbb{C}^m$. Then there exist holomorphic functions $g$ and $h$ such that $hf=g$ on $\mathbb{C}^m$ and $\dim_z h^{-1}(\{0\})\cap g^{-1}(\{0\})\leq m-2$. Therefore the $c$-multiplicity of $f$ is just $\mu^c_f=\mu^0_{g-ch}$ if $c\in\mathbb{C}$ and $\mu^c_f=\mu^0_h$ if $c=\infty$. The function $\mu^c_f: \mathbb{C}^m\to \mathbb{Z}$ is nonnegative and is called the $c$-divisor of $f$.
If $f\not\equiv 0$ on each component of $\mathbb{C}^m$, then $\nu=\mu_f=\mu^0_f-\mu^{\infty}_f$ is called the divisor of $f$.  For $t>0$, the counting function $n_{\nu}$ is defined by
\beas n_{\nu}(t)=t^{-2(m-1)}\int_{A[t]}\nu \upsilon_m^{m-1},\eeas
where $A=\text{supp}\;\nu=\ol{\{z\in G: \nu_f(z)\neq 0\}}$.
The valence function of $\nu$ is defined by
\[N_{\nu}(r)=N_{\nu}(r,r_0)=\int_{r_0}^r n_{\nu}(t)\frac{dt}{t}\;\;(r\geq r_0).\]

\smallskip
We write $N_{\mu_f^a}(r)=N(r,a;f)$ if $a\in\mathbb{C}$ and $N_{\mu_f^a}(r)=N(r,f)$ if $a=\infty$. For $k\in\mathbb{N}$, define the truncated multiplicity functions on $\mathbb{C}^m$ by $\mu_{f,k}^a(z)=\min\{\mu_f^a(z),k\}$, $\mu_{f(k}^a(z)=\mu_f^a(z)$ if $\mu_f^a(z)\geq  k$ and $\mu_{f(k}^a(z)=0$ if $\mu_f^a(z)<k$. Also we define the truncated valence functions $N_{\nu}(t)=\ol N(t,a;f)$ if $\nu=\mu_{f,1}^a$ and $N_{\nu}(t)=N_{(k}(t,a;f)$, if $\nu=\mu_{f(k}^a$.

\smallskip
An algebraic subset $X$ of $\mathbb{C}^m$ is defined as a subset
\[X=\left\lbrace z\in\mathbb{C}^m: P_j(z)=0,\;1\leq j\leq l\right\rbrace\]
with finitely many polynomials $P_1(z),\ldots, P_l(z)$.
A divisor $\nu$ on $\mathbb{C}^m$ is said to be algebraic if $\nu$ is the zero divisor of a polynomial. In this case the counting
function $n_{\nu}$ is bounded (see \cite{GK1,WS1}).

\smallskip
With the help of the positive logarithm function, we define the proximity function of $f$ by
\[m(r, f)=\int_{\mathbb{C}^m\langle r\rangle} \log^+ |f|\;\sigma_m \geq  0.\]

The characteristic function of $f$ is defined by $T(r, f)=m(r,f)+N(r,f)$.
We define $m(r,a;f)=m(r,f)$ if $a=\infty$ and $m(r,a;f)=m(r,1/(f-a))$ if $a\in\mathbb{C}$. Now if $a\in\mathbb{C}$, then the first main theorem becomes $m(r,a;f)+N(r,a;f)=T(r,f) + O(1)$, where $O(1)$ denotes a bounded function when $r$ is sufficiently large.
We define the order of $f$ by
\[\rho(f):=\limsup _{r \rightarrow \infty} \frac{\log T(r, f)}{\log r}.\]

\smallskip
We define the linear measure $m(E):=\int_E dt$, the logarithmic measure $l(E):=\int_{E\cap [1,\infty)} \frac{d t}{t}$ and the upper density measure
\[\ol {\text{dens}}\;E=\lim\limits_{r\rightarrow \infty}\frac{1}{r}\int_{E\cap [1,r]} dt\]
for a set $E\subset [0,\infty)$. Moreover, if $l(E)<+\infty$, resp., $m(E)<+\infty$, then $E$ is of zero upper density.
Let $S(f)=\{g:\mathbb{C}^m\to\mathbb{P}^1\;\text{meromorphic}: \parallel T(r,g)=o(T(r,f))\}$, where $\parallel$ indicates that the equality holds only outside a set $E$ with zero upper density measure and the element in $S(f)$ is called the small function of $f$.

\section{\bf{Key lemmas}}
\begin{lem}\label{L.1} \cite[Lemma 1.37]{HLY1} Let $f$ be a non-constant meromorphic function in $\mathbb{C}^m$ and $I=(\alpha_1,\ldots,\alpha_m)\in \mathbb{Z}^m_+$ be a multi-index. Then for any $\varepsilon>0$, we have
\[m\left(r,\frac{\partial^I(f)}{f}\right)\leq |I|\log^+T(r,f)+|I|(1+\varepsilon)\log^+\log T(r,f)+O(1)\]
for all $r\not\in E$, where $l(E)<+\infty$.
\end{lem}

\begin{lem}\label{L.2} \cite[Lemma 1.2]{HY1} Let $f$ be a non-constant meromorphic function in $\mathbb{C}^m$ and let $a_1,a_2,\ldots,a_q$ be different points in $\mathbb{C}\cup \{\infty\}$. Then
\beas \parallel (q-2)T(r,f)\leq \sideset{}{_{j=1}^{q}}{\sum} \ol N(r,a_j;f)+O(\log (rT(r,f))).\eeas
\end{lem}

\begin{lem}\label{L.3} \cite[Theorem 1.26]{HLY1} Let $f$ be non-constant meromorphic function in $\mathbb{C}^m$. Assume that
$R(z, w)=\frac{A(z, w)}{B(z, w)}$. Then
\beas T\left(r, R_f\right)=\max \{p, q\} T(r, f)+O\Big(\sideset{}{_{j=0}^p}{\sum} T(r, a_j)+\sideset{}{_{j=0}^q}{\sum}T(r, b_j)\Big),\eeas
where $R_f(z)=R(z, f(z))$ and two coprime polynomials $A(z, w)$ and $B(z,w)$ are given
respectively as follows:
\[A(z,w)=\sum_{j=0}^p a_j(z)w^j\;\;\text{and}\;\;B(z,w)=\sum_{j=0}^q b_j(z)w^j.\]
\end{lem}

\begin{lem}\label{L.6}\cite[Lemma 3.2]{HLY1} Let $f_j\not\equiv 0\;(j=1,2,\ldots,n;n\geq 3)$ be meromorphic functions in $\mathbb{C}^m$ such that $f_1,\ldots,f_{n-1}$ are non-constants and $f_1+\cdots+f_n=1$. If
\beas \parallel \;\;\sideset{}{_{j=1}^n}{\sum}\Big\lbrace N_{n-1}(r,0;f_j)+(n-1)\ol{N}(r,f_j)\Big\rbrace<\lambda T(r,f_j)+O(\log^+(T(r,f_j))\eeas
holds for $j=1,2,\ldots,n-1$, where $\lambda<1$, then $f_n\equiv 1$.
\end{lem}

\begin{lem}\label{L.5}\cite[Proposition 3.2]{hy1} Let $P$ be a non-constant entire function in $\mathbb{C}^m$. Then
\[\rho(e^P)=
\begin{cases}
\deg(P) & \text{if $P$ is a polynomial,}\\
+\infty. & \text{otherwise}
\end{cases}\]
\end{lem}

\begin{lem}\label{L.8}\cite[Theorem 2.1]{CX1} Let $f$ be a non-constant meromorphic function in $\mathbb{C}^m$ and let $c\in \mathbb{C}^m$. If
\[\limsup\limits_{r\rightarrow \infty} \frac{\log T(r,f)}{r}=0,\]
then
\beas \parallel\;m\left(r, \frac{f(z+c)}{f(z)}\right)+m\left(r,\frac{f(z)}{f(z+c)}\right)=o(T(r,f)).\eeas
\end{lem}

\begin{lem}\label{L.7}\cite[Theorem 2.2]{CX1} Let $f$ be a non-constant meromorphic function on $\mathbb{C}^m$ with
\[\limsup\limits_{r\rightarrow \infty} \frac{\log T(r,f)}{r}=0,\]
then
\[\parallel\;T(r,f(z+c))=T(r,f)+o(T(r,f))\]
holds for any constant $c\in\mathbb{C}^m$.
\end{lem}

Let $f$ be a non-constant meromorphic function in $\mathbb{C}^m$. Define complex differential-difference polynomials as follows:
\bea\label{cl1} P(f(z))=\sum\limits_{\mathbf{p} \in I} a_{\mathbf{p}}(z) f^{p_0}(z)\big(\partial^{\mathbf{i}_1}(f(z))\big)^{p_1} \cdots\big(\partial^{\mathbf{i}_l}(f(z))\big)^{p_l}f^{p_{l+1}}(z + \hat q_{l+1}) \cdots f^{p_{l+s}}(z + \hat q_{l+s}),\eea
$\mathbf{p}=\left(p_0, \ldots, p_{l+s}\right) \in \mathbb{Z}_{+}^{l+s+1}$,
\bea\label{cl2} Q(f(z))=\sum\limits_{\mathbf{q} \in J} c_{\mathbf{q}}(z) f^{q_0}(z)\big(\partial^{\mathbf{j}_1}(f(z))\big)^{q_1} \cdots\big(\partial^{\mathbf{j}_l}(f(z))\big)^{q_l}f^{q_{l+1}}(z+\tilde q_{l+1}) \cdots f^{q_{l+t}}(z+\tilde q_{l+t}),\eea
$\mathbf{q}=\left(q_0, \ldots, q_{l+t}\right) \in \mathbb{Z}_{+}^{l+t+1}$ and
\bea\label{cl3} B(f(z))=\sideset{}{_{k=0}^{n}}{\sum} b_k(z) f^k(z),\eea
where $I, J$ are finite sets of distinct elements, $a_{\mathbf{p}}, b_k, c_{\mathbf{q}}\in S(f)$ such that $b_n\not\equiv 0$ and $\hat q_i, \tilde q_j\in\mathbb{C}^m$.

For the case, when $P(f(z))$ and $Q(f(z))$ are complex differential polynomials, Hu and Yang \cite[Lemma 2.1]{ps1} generalised Clunie-lemma to high dimension. In 2020, Cao and Xu \cite[Theorem 3.6]{CX1} improved and extended Laine-Yang's difference analogue of Clunie theorem in one variable \cite[Theorem 2.3]{LY1} to high dimension by using Lemma \ref{L.8}. Now by using Lemmas \ref{L.1} and \ref{L.8} and proceeding in the same way as done in the proofs of Lemma 2.1 \cite{ps1} and Theorem 2.3 \cite{LY1} we get the following lemma.

\begin{lem}\label{L.4} Let $f$ be a non-constant meromorphic function in $\mathbb{C}^m$ such that
\[\limsup_{r \to \infty} \frac{\log T(r,f)}{r} = 0.\]

Assume that $f$ satisfies the complex differential-difference equation
\[B(f) Q(f)=P(f),\]
where $P(f)$, $Q(f)$ and $B(f)$ are defined as in (\ref{cl1}), (\ref{cl2}) and (\ref{cl3}) respectively.
If $\deg(P(f)) \leq n=\deg(B(f))$, then
\[\parallel\;m(r, Q(f))=o(T(r, f)).\]
\end{lem}

\begin{lem}\label{L.8}\cite[Lemma 1.68]{HLY1} Let $f_1:\mathbb{C}^m\to\mathbb{P}^1$ and $f_2:\mathbb{C}^m\to\mathbb{P}^1$ be two non-constant meromorphic functions. Then for $r>0$ we have
\[N(r,0;f_1f_2)-N(r,f_1f_2)=N(r,0;f_1)+N(r,0;f_2)-N(r,f_1)-N(r,f_2).\]
\end{lem}

\section {{\bf Proof of Theorem \ref{t2.1}}}
\begin{proof} Let $(f_1(z),f_2(z))$ be a pair of finite order non-constant entire solution of the system of equations (\ref{pds1}). Set
\bea\label{lm} F_1(z)=\frac{\partial f_1(z)}{\partial z_1}\;\;\text{and}\;\;F_2(z)=\frac{\partial f_2(z)}{\partial z_1}.\eea

Since $f_i(z)$ are entire functions for $i=1,2$, applying Lemma \ref{L.1} to (\ref{lm}), we have
\bea\label{lm1} T(r,F_i)=m(r,F_i)\leq m\left(r,\frac{F_i}{f_i}\right)+m(r,f_i)\leq T(r,f_i)+o(T(r,f_i))\eea
and so $o(T(r,F_i))$ can be replaced by $o(T(r,f_i))$, where $i=1,2$. Now using Lemma \ref{L.3} to the system of equations (\ref{pds1}), we get
\bea\label{lm.0} \parallel\;n_i T(r,F_i)+o(T(r,F_i))=m_i T(r,f_j(z+c))+o(T(r,f_j(z+c))),\eea
where $i,j\in\{1,2\}$ such that $i\neq j$. Since $\rho(f_i)<+\infty$, we have
\[\limsup_{r \to \infty} \frac{\log T(r,f_i)}{r}=0\]
for $i=1,2$ and so by Lemma \ref{L.7}, we have
\bea\label{lm.00} \parallel\;T(r,f_i(z+c))=T(r,f_i)+o(T(r,f_i)),\eea
where $i\in\{1,2\}$. Now from (\ref{lm.0}) and (\ref{lm.00}), we get
\bea\label{lm.1} \parallel\;n_i T(r,F_i)+o(T(r,F_i))=m_i T(r,f_j)+o(T(r,f_j)),\eea
where $i,j\in\{1,2\}$ such that $i\neq j$. Clearly from (\ref{lm1}) and (\ref{lm.1}), we obtain
\bea\label{lm2} m_i T(r,f_j)+o(T(r,f_j))\leq n_i T(r,f_i)+o(T(r,f_i)),\eea
where $i,j\in\{1,2\}$ such that $i\neq j$. Consequently (\ref{lm2}) shows that $o(T(r,f_i))$ can be replaced by $o(T(r,f_j))$, where $i,j\in\{1,2\}$ such that $i\neq j$. On the other hand, using Lemmas \ref{L.1}, \ref{L.3}, (\ref{lm.00}) and (\ref{lm.1}) to (\ref{pds1}), we obtain
\bea\label{lm.2}\parallel\;m_i T(r, f_j))&=& m_i T(r, f_j(z+c))+o(T(r, f_j))\\
& =&T\left(r, f_j^{m_i}(z+c)\right)+o(T(r, f_j))\nonumber\\
& =&T\left(r,1-F_i^{n_i}\right)+o(T(r, f_j)) \nonumber\\
& =&n_i T(r,F_i)+o(T(r,F_i))+o(T(r, f_j)) \nonumber\\
& =&n_i m(r,F_i)+o(T(r, F_i))+o(T(r, f_j)) \nonumber\\
& \leq & n_i\left(m\left(r, \frac{F_i}{f_i}\right)+m(r, f_i)\right)+o(T(r, F_i))+o(T(r, f_j))\nonumber \\
& =&n_i T(r, f_i)+o(T(r, f_j)),\nonumber\eea
where $i,j\in\{1,2\}$ such that $i\neq j$.\par

\smallskip
First we suppose that $m_1m_2-n_1n_2>0$. Then from  (\ref{lm.2}), one can easily deduce that
\[\parallel\;(m_1m_2-n_1n_2)T(r,f_j)\leq o(T(r, f_j)),\]
which is impossible, where $j\in\{1,2\}$. Therefore in this case, the system of Fermat-type partial differential-difference equations (\ref{pds1}) does not have any finite order transcendental entire solutions.\par

\medskip
Next we suppose that $m_1m_2-n_1n_2\leq 0$.
Let
\bea\label{lm.3} h_i(z)=\frac{F_i^{n_i}(z)-1}{F_i^{n_i}(z)},\eea
where $i\in\{1,2\}$. Clearly $h_i(z)$ is a non-constant meromorphic function in $\mathbb{C}^m$. Again using Lemma \ref{L.3} to (\ref{lm.3}), we get
\bea\label{lm.4}\parallel T(r,h_i)+o(T(r,h_i))=n_iT(r,F_i)+o(T(r,F_i)),\eea
where $i\in\{1,2\}$. Now using (\ref{pds1}) to (\ref{lm.3}), one can easily deduce that
\[\parallel\;\ol N(r,h_i)\leq \ol N(r,0, F_i^{n_i})=\ol N(r,0,F_i)+o(T(r,F_i)),\]
\[\parallel\;\ol N(r,0,h_i)=\ol N(r,1,F_i^{n_i})\leq \ol N(r,0,f_j^{m_i}(z+c))=\ol N(r,0,f_j(z+c))+o(T(r,f_j(z+c)))\]
and $\parallel \ol N(r,1,h_i)=0$,
where $i,j\in\{1,2\}$ such that $i\neq j$. Therefore in view of the first main theorem and using Lemma \ref{L.2} and (\ref{lm.4}), we get
\bea\label{lm.6} \parallel\;n_i T(r,F_i)&=&T(r,h_i)+o(T(r,h_i))\\&\leq& \ol N(r,h_i)+\ol N(r,0,h_i)+\ol N(r,1,h_i)+o(T(r,h_i))\nonumber\\&\leq&
\ol N(r,0,F_i)+\ol N(r,0,f_j(z+c))+o(T(r,F_i))+o(T(r,f_j(z+c)))\nonumber\\&\leq&
T(r,F_i)+T(r,f_j(z+c))+o(T(r,F_i))+o(T(r,f_j(z+c))),\nonumber
\eea
where $i,j\in\{1,2\}$ such that $i\neq j$. Using (\ref{lm.1}) to (\ref{lm.6}), we get
\[\parallel\;\left(n_i-1-\frac{n_i}{m_i}\right)T(r,F_i)\leq o(T(r,F_i)),\]
where $i\in\{1,2\}$ and so we have
\bea\label{lm.5}\frac{1}{n_i}+\frac{1}{m_i}\geq 1,\eea
where $i\in\{1,2\}$.
We consider the following three cases.\par

\medskip
{\bf Case 1.} Let $n_1= m_1$. Note that $n_1+m_1>2$ and so (\ref{lm.5}) gives $n_1=m_1=2$. Since $m_1m_2-n_1n_2\leq 0$, it follows that $m_2\leq n_2$. Now we consider following two sub-cases.\par

\medskip
{\bf Sub-case 1.1.} Let $m_2<n_2$. In this case from (\ref{lm.5}), we obtain $m_2=1$ and $n_2>1$. Since $m_1=n_1=2$, from (\ref{pds1}),  we get
\bea\label{lm.7} F_1^2(z)+f_2^2(z+c)=1\eea
and
\bea\label{lm.8} F_2^{n_2}(z)+f_1(z+c)=1.\eea

Clearly (\ref{lm.8}) gives
\beas n_2F_2^{n_2-1}(z)\frac{\partial F_2(z)}{\partial z_1}+F_1(z+c)=0\eeas
and so from (\ref{lm.7}), we get
\bea\label{lm.9}\left(n_2F_2^{n_2-1}(z)\frac{\partial F_2(z)}{\partial z_1}\right)^2+f_2^2(z+2c)=1.\eea

Now using Theorem 1.2.A to (\ref{lm.9}), one has
\bea\label{lm.10} n_2F_2^{n_2-1}(z)\frac{\partial F_2(z)}{\partial z_1}=\cos (h(z))\eea
and
\bea\label{lm.11} f_2(z+2c)=\sin (h(z)),\eea
where $h:\mathbb{C}^m\to \mathbb{C}$ is a non-constant entire function. Since $\rho(f_i)<+\infty$ for $i=1,2$, using Lemmas \ref{L.3}, \ref{L.5} and (\ref{lm.00}) to (\ref{lm.11}), we conclude that $h(z)$ is a non-constant polynomial in $\mathbb{C}^m$. Clearly from (\ref{lm.11}), we obtain
\bea\label{lm.12}
\begin{cases}
F_2(z+2c)=\frac{\partial h(z)}{\partial z_1}\cos (h(z)),\\
\frac{\partial F_2(z+2c)}{\partial z_1}=\frac{\partial^2 h(z)}{\partial z_1^2}\cos (h(z))-\left(\frac{\partial h(z)}{\partial z_1}\right)^2\sin (h(z)).
\end{cases}
\eea

Therefore from (\ref{lm.10}) and (\ref{lm.12}), we have
\bea &&\label{lm.13}n_2\left(\left(\frac{\partial h(z)}{\partial z_1}\right)^{n_2-1}\frac{\partial^2 h(z)}{\partial z_1^2}\cos (h(z))-\left(\frac{\partial h(z)}{\partial z_1}\right)^{n_2+1}\sin (h(z))\right)\cos^{n_2-1} (h(z))\\&=&\cos (h(z+2c)).\nonumber\eea

Clearly $\frac{\partial h(z)}{\partial z_1}\not\equiv 0$. Note that the sets of multiple zeros of both $\cos (h(z))$ and $\cos (h(z+2c))$ are algebraic. But if $n_2\geq 3$, then the set of multiple zeros of $\cos^{n_2-1} (h(z))$ is not algebraic and so from (\ref{lm.13}), we get a contradiction. Hence $n_2=2$ and so (\ref{lm.13}) yields
\beas 2\left(\frac{\partial h(z)}{\partial z_1}\frac{\partial^2 h(z)}{\partial z_1^2}\cos (h(z))-\left(\frac{\partial h(z)}{\partial z_1}\right)^{3}\sin (h(z))\right)\cos (h(z))=\cos (h(z+2c)),\eeas
i.e.,
\bea\label{lm.14} 2\frac{\partial h(z)}{\partial z_1}\frac{\partial^2 h(z)}{\partial z_1^2}\cos^2(h(z))- \bigg(\frac{\partial h(z)}{\partial z_1}\bigg)^3\sin (2h(z))=\cos(h(z+2c)).\eea

\medskip
First we suppose that $\frac{\partial^2 h(z)}{\partial z_1^2}\equiv 0$. Clearly $\frac{\partial h(z)}{\partial z_1}=d(z_2,z_3,\cdots,z_m)$, a polynomial in $\mathbb{C}^{m-1}$. Then from (\ref{lm.14}), we get
\bea\label{lm.14a}\iota d^3e^{\iota (2h(z)+h(z+2c))}-\iota d^3e^{\iota (-2h(z)+h(z+2c))}-e^{2\iota h(z+2c)}=1.\eea

Since both $2h(z)+h(z+2c)$ and $-2h(z)+h(z+2c)$ are non-constants, applying Lemma \ref{L.6} to (\ref{lm.14a}), we get a contradiction.\par

\medskip
Next we suppose that $\frac{\partial^2 h(z)}{\partial z_1^2}\not\equiv 0$. By simple calculations on (\ref{lm.14}), we get
\bea\label{lm.15} -\frac{1}{\frac{\partial h(z)}{\partial z_1}\frac{\partial^2 h(z)}{\partial z_1^2}}e^{\iota h(z+2c)}-\frac{1}{\frac{\partial h(z)}{\partial z_1}\frac{\partial^2 h(z)}{\partial z_1^2}}e^{-\iota h(z+2c)}+\frac{1}{2}G_1(z)e^{2\iota h(z)}+\frac{1}{2}G_2(z)e^{-2\iota h(z)}=1,\eea
where
\[G_1(z)=1-2\iota\frac{\big(\frac{\partial h(z)}{\partial z_1}\big)^2}{\frac{\partial^2 h(z)}{\partial z_1^2}}\;\;\text{and}\;\; G_2(z)=1+2\iota\frac{\big(\frac{\partial h(z)}{\partial z_1}\big)^2}{ \frac{\partial^2 h(z)}{\partial z_1^2}}.\]

If $G_1(z)\equiv 0$, then from (\ref{lm.15}), we get
\[-\frac{1}{\frac{\partial h(z)}{\partial z_1}\frac{\partial^2 h(z)}{\partial z_1^2}}e^{\iota h(z+2c)}-\frac{1}{\frac{\partial h(z)}{\partial z_1}\frac{\partial^2 h(z)}{\partial z_1^2}}e^{-\iota h(z+2c)}+e^{-2\iota h(z)}=1\]
and so using Lemma \ref{L.6}, we get a contradiction. Hence $G_1(z)\not\equiv 0$. Similarly we can prove that $G_2(z)\not\equiv 0$.
Finally using Lemma \ref{L.6} to (\ref{lm.15}), we get a contradiction.\par

\medskip
{\bf Sub-case 1.2.} Let $n_2=m_2$. Since $n_2+m_2>2$, from (\ref{lm.5}), we get $n_2=m_2=2$. Consequently from (\ref{pds1}), we obtain
$F_1^2(z)+f_2^2(z+c)=1$ and $F_2^2(z)+f_1^2(z+c)=1$, i.e.,
\bea\label{lm.16} (F_1(z)+\iota f_2(z+c))(F_1(z)-\iota f_2(z+c))=1\eea
and
\bea\label{lm.17} (F_2(z)+\iota f_1(z+c))(F_2(z)-\iota f_1(z+c))=1.\eea

Clearly from (\ref{lm.16}) and (\ref{lm.17}), we see that the functions both $F_1(z)+\iota f_2(z+c)$, $F_1(z)-\iota f_2(z+c)$, $F_2(z)+\iota f_1(z+c)$ and $F_2(z)-\iota f_1(z+c)$ have no zeros. So we may assume that
\bea\label{sm.3} F_1(z)+\iota f_2(z+c)=e^{\iota P_1(z)},\eea
\bea\label{sm.4} F_1(z)-\iota f_2(z+c)=e^{-\iota P_1(z)},\eea
\bea\label{dpm4} F_2(z)+\iota f_1(z+c)=e^{\iota P_2(z)}\eea
and
\bea\label{dpm5} F_2(z)-\iota f_1(z+c)=e^{-\iota P_2(z)},\eea
where $P_1(z)$ and $P_2(z)$ are  polynomials in $\mathbb{C}^m$. Solving (\ref{sm.3}) and (\ref{sm.4}), we get
\bea\label{sm.5} F_1(z)=\frac{\partial f_1(z)}{\partial z_1}=\frac{e^{\iota P_1(z)}+e^{-\iota P_1(z)}}{2}=\cos (P_1(z))\eea
and
\bea\label{sm.6} f_2(z+c)=\frac{e^{\iota P_1(z)}-e^{-\iota P_1(z)}}{2\iota}=\sin (P_1(z)).\eea

Again solving (\ref{dpm4}) and (\ref{dpm5}), we obtain
\bea\label{dpm6} F_2(z)=\frac{\partial f_2(z)}{\partial z_1}=\frac{e^{\iota P_2(z)}+e^{-\iota P_2(z)}}{2}=\cos (P_2(z))\eea
and
\bea\label{dpm7} f_1(z+c)=\frac{e^{\iota P_2(z)}-e^{-\iota P_2(z)}}{2\iota}=\sin (P_2(z)).\eea

Clearly from (\ref{sm.6}) and (\ref{dpm6}), we have
\bea\label{sm.8} \frac{\partial P_1(z)}{\partial z_1}e^{\iota (P_1(z)+P_2(z+c))}+\frac{\partial P_1(z)}{\partial z_1}e^{-\iota (P_1(z)-P_2(z+c))}-e^{2\iota P_2(z+c)}=1,\eea
which implies that $\frac{\partial P_1(z)}{\partial z_1}\not\equiv 0$. Again from (\ref{sm.5}) and (\ref{dpm7}), we get
\bea\label{dp9} \frac{\partial P_2(z)}{\partial z_1}e^{\iota (P_2(z)+P_1(z+c))}+\frac{\partial P_2(z)}{\partial z_1}e^{-\iota (P_2(z)-P_1(z+c))}-e^{2\iota P_1(z+c)}=1,\eea
which shows that $\frac{\partial P_2(z)}{\partial z_1}\not\equiv 0$.
Thus, using Lemma \ref{L.6} to (\ref{sm.8}) and (\ref{dp9}), we have respectively
\beas \text{either}\;\;\frac{\partial P_1(z)}{\partial z_1}e^{\iota (P_1(z)+P_2(z+c))}=1 \;\;\text{or}\;\; \frac{\partial P_1(z)}{\partial z_1}e^{-\iota (P_1(z)-P_2(z+c))}=1\eeas
and
\beas \text{either}\;\;\frac{\partial P_2(z)}{\partial z_1}e^{\iota (P_2(z)+P_1(z+c))}=1 \;\;\text{or}\;\; \frac{\partial P_2(z)}{\partial z_1}e^{-\iota (P_2(z)-P_1(z+c))}=1.\eeas

Next, we consider the following four sub-cases.\par

\medskip
{\bf Sub-case 1.2.1.} Let
\bea\label{dpm10} \begin{cases}  \frac{\partial P_1(z)}{\partial z_1}e^{\iota (P_1(z)+P_2(z+c))}=1\\ \frac{\partial P_2(z)}{\partial z_1}e^{\iota (P_2(z)+P_1(z+c))}=1.\end{cases}\eea

It follows from (\ref{dpm10}) that $P_1(z+c)+P_2(z)$ and $P_2(z+c)+P_1(z)$ are both constants. Consequently both $P_1(z+2c)-P_1(z)$ and $P_2(z+2c) -P_2(z)$ are constants. Assume that
\bea\label{1.3.1} e^{-\iota (P_1(z)+P_2(z+c))}=A\;\;\text{and}\;\;e^{-\iota (P_2(z)+P_1(z+c))}=B.\eea

Then from (\ref{dpm10}), we have
\[\frac{\partial P_1(z)}{\partial z_1}=A\;\;\text{and}\;\;\frac{\partial P_2(z)}{\partial z_1}=B,\]
where $A^2=1$ and $B^2=1$. Obviously $P_1(z)=Az_1+\tilde Q_1(z_2,\ldots,z_m)$ and $P_2(z)=Bz_1+\tilde Q_2(z_2,\ldots,z_m)$,
where $\tilde Q_1(z_2,\ldots,z_m)$ and $\tilde Q_2(z_2,\ldots,z_m)$ are polynomials in $\mathbb{C}^{m-1}$. Since both $P_1(z+2c)-P_1(z)$ and $P_2(z+2c) -P_2(z)$ are constants, it follows that $\tilde Q_i(z_2+2c_2,\ldots,z_m+2c_m)-Q_i\tilde (z_2,\ldots,z_m)$ is also a constant for $i=1,2$. Therefore we may assume that
\bea\label{dpm11}
\begin{cases}
P_1(z)=Az_1+A_{12}z_2+\ldots+A_{1m}z_m+\hat Q_1(z_2,\ldots,z_m), \\
P_2(z)=Bz_1+B_{12}z_2+\ldots+B_{1m}z_m+\hat Q_2(z_2,\ldots,z_m),
\end{cases}
\eea
where $A_{1i}, B_{1i}\in\mathbb{C}$ for $i=2,\ldots,m$ and $\hat Q_i(z_2,\ldots,z_m)$ is a polynomial in $\mathbb{C}^{m-1}$ such that
$\hat Q_i(z_2+2c_2,\ldots,z_m+2c_m)=\hat Q_i(z_2,\ldots,z_m)$
for $i=1,2$. On the other hand, (\ref{1.3.1}) gives
\[\frac{B}{A}e^{\iota (P_1(z+2c)-P_1(z))}=1\;\;\text{and}\;\;\frac{A}{B}e^{\iota (P_2(z+2c)-P_2(z))}=1.\]

Clearly from (\ref{dpm11}), we have respectively
\[\frac{B}{A}e^{2\iota (Ac_1+A_{12}c_2+\ldots+A_{1m}c_m)}=1\;\;\text{and}\;\;\frac{A}{B}e^{2\iota (Bc_1+B_{12}c_2+\ldots+B_{1m}c_m)}=1.\]

Finally, from (\ref{sm.6}) and (\ref{dpm7}), we may assume that
\[(f_1(z),f_2(z))=\sin (L_1(z)+Q_1(z_2,z_3,\ldots, z_m)),\sin (L_2(z)+Q_2(z_2,z_3,\ldots, z_m))),\]
where $L_1(z)=Az_1+A_{12}z_2+\ldots+A_{1m}z_m$, $L_2(z)=Bz_1+B_{12}z_2+\ldots+B_{1m}z_m$, $A_{1i},B_{1i}\in\mathbb{C}$ for $i=2,\ldots,m$ such that $A^2=1$, $B^2=1$, $\frac{B}{A}e^{2\iota (Ac_1+A_{12}c_2+\ldots+A_{1m}c_m)}=1$ and $\frac{A}{B}e^{2\iota (Bc_1+B_{12}c_2+\ldots+B_{1m}c_m)}=1$ and
$Q_i(z_2,\ldots,z_m)$ is a polynomial such that for $i=1,2$
\[Q_i(z_2+2c_2,\ldots,z_m+2c_m)=Q_i(z_2,\ldots,z_m).\]

\medskip
{\bf Sub-case 1.2.2.} Let
\bea\label{dpm14} \begin{cases}  \frac{\partial P_1(z)}{\partial z_1}e^{\iota (P_1(z)+P_2(z+c))}=1\\ \frac{\partial P_2(z)}{\partial z_1}e^{-\iota (P_2(z)-P_1(z+c))}=1.\end{cases}\eea

Clearly from (\ref{dpm14}), it follows that $P_1(z+c)+P_2(z)$ and $P_2(z+c)-P_1(z)$ are both constants. This means $P_2(z+2c)+P_2(z)$ is also a constant, which contradicts the fact that $P_2(z)$ is a non-constant polynomial.\par

\medskip
{\bf Sub-case 1.2.3.} Let
 \bea\label{dpm15} \begin{cases} \frac{\partial P_1(z)}{\partial z_1}e^{-\iota (P_1(z)-P_2(z+c))}=1\\ \frac{\partial P_2(z)}{\partial z_1}e^{\iota (P_2(z)+P_1(z+c))}=1.\end{cases}\eea

We deduce from (\ref{dpm15}) that $P_1(z)-P_2(z+c)$ and $P_2(z)+P_1(z+c)$ are both constants. This means $P_1(z+2c)+P_1(z)$ is also a constant, which is a contradiction.\par

\medskip
{\bf Sub-case 1.2.4.} Let
\bea\label{dpm16} \begin{cases}  \frac{\partial P_1(z)}{\partial z_1}e^{-\iota (P_1(z)-P_2(z+c))}=1\\ \frac{\partial P_2(z)}{\partial z_1}e^{-\iota (P_2(z)-P_1(z+c))}=1.\end{cases}\eea

Clearly from (\ref{dpm16}), it follows that $P_1(z+c)-P_2(z)$ and $P_2(z+c)-P_1(z)$ are both constants. Consequently $P_1(z+2c)-P_1(z)$ and $P_2(z+2c) -P_2(z)$ are also constants. Assume that
\bea\label{a1.3.1} e^{\iota (P_1(z)-P_2(z+c))}=\tilde A\;\;\text{and}\;\;e^{\iota (P_2(z)-P_1(z+c))}=\tilde B.\eea

Then from (\ref{dpm16}), we have
\[\frac{\partial P_1(z)}{\partial z_1}=\tilde A\;\;\text{and}\;\;\frac{\partial P_2(z)}{\partial z_1}=\tilde B,\]
where $\tilde A^2=1$ and $\tilde B^2=1$. Clearly from (\ref{a1.3.1}), we have $\tilde A \tilde B e^{\iota (P_i(z+2c)-P_i(z))}=1$ for $i=1,2$. Now proceeding in the same way as done in the proof of Sub-case 1.3.1, we can conclude that
\[(f_1(z),f_2(z))=\sin (\tilde L_1(z)+R_1(z_2,z_3,\ldots, z_m)),\sin (\tilde L_2(z)+R_2(z_2,z_3,\ldots, z_m))),\]
where $\tilde L_1(z)=\tilde Az_1+\tilde A_{12}z_2+\ldots+\tilde A_{1m}z_m$, $\tilde L_2(z)=\tilde Bz_1+\tilde B_{12}z_2+\ldots+\tilde B_{1m}z_m$, $\tilde A_{1i}, \tilde B_{1i}\in\mathbb{C}$, $i=2,\ldots,m$ such that $\tilde A^2=1$, $\tilde B^2=1$, $\tilde A \tilde Be^{2\iota (\tilde Ac_1+\tilde A_{12}c_2+\ldots+\tilde A_{1m}c_m)}=1$, $\tilde A \tilde Be^{2\iota (\tilde Bc_1+\tilde B_{12}c_2+\ldots+\tilde B_{1m}c_m)}=1$ and $R_i(z_2,\ldots,z_m)$ is a polynomial such that for $i=1,2$
\[R_i(z_2+2c_2,\ldots,z_m+2c_m)=R_i(z_2,\ldots,z_m).\]

\medskip
{\bf Case 2.} Let $m_1>n_1$. In this case from (\ref{lm.5}), we obtain $m_1>1$ and $n_1=1$. Since $m_1m_2-n_1n_2\leq 0$, it follows that $m_1m_2-n_2\leq 0$. We now consider following sub-cases.

\medskip
{\bf Sub-case 2.1.} Let $m_2>n_2$. Then (\ref{lm.5}) implies $m_2>1$ and $n_2=1$. Therefore $m_1m_2-1\leq 0$, i.e., $m_1m_2\leq 1$. Since $m_1+n_1>2$ and $m_2+n_2>2$, we get a contradiction.\par

\medskip
{\bf Sub-case 2.2.} Let $n_2>m_2$ such that $n_2>2$ and $m_1>\frac{n_2}{n_2-2}$. Clearly (\ref{lm.5}) implies $n_2>1$ and $m_2=1$. In this case, from (\ref{pds1}), we get
\bea\label{llm.1} F_1(z)+f_2^{m_1}(z+c)=1\eea
and
\bea\label{llm.2} F_2^{n_2}(z)+f_1(z+c)=1.\eea

By a simple calculation, (\ref{llm.1}) and (\ref{llm.2}) give
\bea\label{llm.3}-n_2F_2^{n_2-1}(z)\frac{\partial F_2(z)}{\partial z_1}+f_2^{m_1}(z+2c)=1.\eea

Let $\phi(z)=-n_2F_2^{n_2-1}(z)\frac{\partial F_2(z)}{\partial z_1}$. Now using Lemmas \ref{L.1}, \ref{L.3} and \ref{L.8}, we have
\beas T(r,F_2^{n_2}(z))=m(r,F_2^{n_2}(z))&=&m\left(r,\phi(z)\frac{F_2(z)}{\frac{\partial F_2(z)}{\partial z_1}}\right)\\&\leq&
m(r,\phi(z))+m\left(r,\frac{F_2(z)}{\frac{\partial F_2(z)}{\partial z_1}}\right)\\&\leq&
T(r,\phi(z))+T\left(r,\frac{F_2(z)}{\frac{\partial F_2(z)}{\partial z_1}}\right)-N\left(r,\frac{F_2(z)}{\frac{\partial F_2(z)}{\partial z_1}}\right)+O(1)\\&\leq&
T(r,\phi(z))+T\left(r,\frac{\frac{\partial F_2(z)}{\partial z_1}}{F_2(z)}\right)-N\left(r,\frac{F_2(z)}{\frac{\partial F_2(z)}{\partial z_1}}\right)+O(1)\\&\leq&
T(r,\phi)+N\left(r,\frac{\frac{\partial F_2(z)}{\partial z_1}}{F_2(z)}\right)-N\left(r,\frac{F_2(z)}{\frac{\partial F_2(z)}{\partial z_1}}\right)+o(T(r,F_2))\\&\leq&
T(r,\phi)+N(r,0;F_2)-N\left(r,0;\frac{\partial F_2(z)}{\partial z_1}\right)+o(T(r,F_2)),\nonumber
\eeas
i.e.,
\bea\label{llm.4} (n_2-1)T(r,F_2)\leq T(r,\phi)-N\left(r,0;\frac{\partial F_2(z)}{\partial z_1}\right)+o(T(r,F_2)).\eea

Therefore using Lemma \ref{L.2}, (\ref{llm.3}) to (\ref{llm.4}), we obtain
\bea\label{llm.5} (n_2-1)T(r,F_2)&\leq& \ol N(r,0;\phi)+\ol N(r,1;\phi)-N\left(r,0;\frac{\partial F_2(z)}{\partial z_1}\right)+o(T(r,F_2))\\&\leq&\ol N(r,0;F_2)+\ol N(r,0;f_2(z+2c))+o(T(r,F_2))\nonumber\\&\leq&
T(r,F_2)+T(r,f_2(z+2c))+o(T(r,F_2)).\nonumber
\eea

On the other hand, applying Lemmas \ref{L.1} and \ref{L.3} to (\ref{llm.3}), we get
\bea\label{llm.6} m_1 T(r,f_2(z+2c))=T(r,\phi)=m(r,\phi)&\leq& n_2 m(r,F_2)+o(T(r,F_2))\\&=&n_2 T(r,F_2)+o(T(r,F_2)).\nonumber\eea

Consequently (\ref{llm.5}) and (\ref{llm.6}) yield
\[\left(n_2-2-\frac{n_2}{m_1}\right)T(r,F_2)\leq o(T(r,F_2)),\]
which contradicts the fact that $m_1>\frac{n_2}{n_2-2}$.\par

\medskip
{\bf Sub-case 2.3.} Let $n_2=m_2$. Then (\ref{lm.5}) gives $n_2=m_2=2$ and so from (\ref{pds1}), we get
\bea\label{llm.7} F_2^2(z)+f_1^2(z+c)=1.\eea

Now using Theorem 1.2.A to (\ref{llm.7}), one has
\bea\label{llm.8} F_2(z)=\cos (h(z))\eea
and
\bea\label{llm.9} f_1(z+c)=\sin (h(z)),\eea
where $h:\mathbb{C}^m\to \mathbb{C}$ is a non-constant polynomial. Clearly from (\ref{llm.8}) and (\ref{llm.9}), we obtain
\bea\label{llm.10}
\begin{cases}
F_1(z+c)=\frac{\partial h(z)}{\partial z_1}\cos (h(z))=\frac{\partial h(z)}{\partial z_1}F_2(z),\\
\frac{\partial F_1(z+c)}{\partial z_1}=\frac{\partial^2 h(z)}{\partial z_1^2}F_2(z)+\frac{\partial h(z)}{\partial z_1}\frac{\partial F_2(z)}{\partial z_1}.
\end{cases}
\eea

Clearly $\frac{\partial h(z)}{\partial z_1}\not\equiv 0$. On the other hand, (\ref{llm.1}) yields
\beas \frac{\partial F_1(z)}{\partial z_1}+m_1f_2^{m_1-1}(z+c)F_2(z+c)=0,\eeas
i.e.,
\bea\label{llm.11} \left(\frac{\partial F_1(z)}{\partial z_1}\right)^{m_1}=(-m_1)^{m_1}(F_1(z)-1)^{m_1-1}F_2^{m_1}(z+c).\eea

Now using (\ref{llm.10}) to (\ref{llm.11}), we get
\bea\label{llm.12}\left(\frac{\partial^2 h(z)}{\partial z_1^2}F_2(z)+\frac{\partial h(z)}{\partial z_1}\frac{\partial F_2(z)}{\partial z_1}\right)^{m_1}=(-m_1)^{m_1}\left(\frac{\partial h(z)}{\partial z_1}F_2(z)-1\right)^{m_1-1}F_2^{m_1}(z+2c).\eea

Applying Lemma \ref{L.4} to (\ref{llm.12}), we deduce that
\[m\left(r,\left(\frac{\partial h(z)}{\partial z_1}F_2(z)-1\right)^{m_1-1}\right)=o(T(r,F_2))\]
and so by Lemma \ref{L.3}, we get $(m_1-1)T(r,F_2)=o(T(r,F_2))$, which is impossible.\par

\medskip
{\bf Case 3.} Let $n_1>m_1$. In this case from (\ref{lm.5}), we obtain $m_1=1$ and $n_1>1$. Since $m_1m_2-n_1n_2\leq 0$, it follows that $m_2-n_1n_2\leq 0$.
We now consider following sub-cases.

\medskip
{\bf Sub-case 3.1.} Let $m_2>n_2$ such that $n_1\geq 3$ and $m_2>\frac{n_1}{n_1-2}$. Therefore (\ref{lm.5}) implies that $m_2>1$ and $n_2=1$. In this case, from (\ref{pds1}), we get
\bea\label{lllm.1} F_1^{n_1}(z)+f_2(z+c)=1\eea
and
\bea\label{lllm.2} F_2(z)+f_1^{m_2}(z+c)=1.\eea

By simple calculations on (\ref{lllm.1}) and (\ref{lllm.2}) give
\beas-n_1F_1^{n_1-1}(z)\frac{\partial F_1(z)}{\partial z_1}+f_1^{m_2}(z+2c)=1.\eeas

Since $m_2>\frac{n_1}{n_1-2}$, proceeding in the same way as done in the proof of Sub-case 2.2, we get a contradiction.\par

\medskip
{\bf Sub-case 3.2.} Let $m_2<n_2$. In this case, (\ref{lm.5}) implies that $m_2=1$ and $n_2>1$. Since $m_2-n_1n_2\leq 0$, it follows that $n_1n_1\ge 1$. Therefore (\ref{pds1}) gives
\bea\label{md1}F_i^{n_i}(z)+f_j(z+c)=1\eea
for $i,j\in\{1,2\}$ such that $i\neq j$.
Differentiating (\ref{md1}) partially with respect to $z_1$, we get
\bea\label{md2} n_i F_i^{n_i-1}(z)\frac{\partial F_i(z)}{\partial z_1}+F_j(z+c)=0.\eea
for $i,j\in\{1,2\}$ such that $i\neq j$.
Clearly (\ref{md1}) and (\ref{md2}), we obtain
\bea\label{md3}\left(-n_iF_i^{n_i-1}(z)\frac{\partial F_i(z)}{\partial z_1}\right)^{n_j}+f_i(z+2c)=1.\eea
for $i,j\in\{1,2\}$ such that $i\neq j$.
Differentiating (\ref{md3}) partially with respect to $z_1$, we get
\bea\label{md4} F_i^{n_1n_2-n_j-1}(z)\alpha_i(z) =-F_i(z+2c),\eea
where
\bea\label{md5}\alpha_i(z)= (-n_i)^{n_j}n_j\left(\frac{\partial F_i(z)}{\partial z_1}\right)^{n_j-1}\left[(n_i-1)\left(\frac{\partial F_i(z)}{\partial z_1}\right)^2+F_i(z)\frac{\partial^2 F_i(z)}{\partial z_1^2}\right].\eea
for $i,j\in\{1,2\}$ such that $i\neq j$.
Clearly $\alpha_i\not\equiv 0$, otherwise from (\ref{md5}), we get $F_i(z+2c)\equiv 0$, i.e., $F_i(z)\equiv 0$, which is impossible for $i=1,2$. Since $n_i>1$, we have $n_1n_2-n_i\geq 2$ for $i=1,2$. Therefore applying Lemma \ref{L.4} to (\ref{md4}), we get $\parallel m(r,\alpha_i)=o(T(r,F_i))$ and so $\parallel T(r,\alpha_i)=o(T(r,F_i))$ for $i=1,2$.
Now we have to consider following sub-cases.\par

\medskip
{\bf Sub-case 3.2.1.} Let $n_1n_2-n_i>2$ for $i=1,2$. Now using Lemma \ref{L.4} to (\ref{md4}), we get $\parallel m(r,F_i\alpha)=o(T(r,F_i))$ and so $\parallel T(r,F_i\alpha)=o(T(r,F_i))$ for $i=1,2$. Since $\parallel T(r,\alpha_i)=o(T(r,F_i))$, it follows that
\[\parallel\;T(r,F_i)\leq T(r,F_i\alpha_i)+T\left(r,\frac{1}{\alpha_i}\right)=T(r,F_i\alpha_i)+T(r,\alpha_i)=o(T(r,F_i)),\]
which is impossible for $i=1,2$.\par

\medskip
{\bf Sub-case 3.2.2.} Let $n_1n_2-n_i=2$ for $i=1,2$. In this case, we must have $n_1=n_2=2$. Then from (\ref{md4}) and (\ref{md5}), we have respectively
\bea\label{md7} -F_i(z)\alpha_i(z)=F_i(z+2c)\eea
for $i=1,2$ and
\bea\label{md8a}\alpha_i(z)=8\frac{\partial F_i(z)}{\partial z_1}\left[\left(\frac{\partial F_i(z)}{\partial z_1}\right)^2+F_i(z)\frac{\partial^2 F_i(z)}{\partial z_1^2}\right]=8\frac{\partial F_i(z)}{\partial z_1}\frac{\partial \left(F_i(z)\frac{\partial F_i(z)}{\partial z_1}\right)}{\partial z_1},\eea
i.e.,
\bea\label{md8}\alpha_i(z)F_i(z)=8F_i(z)\frac{\partial F_i(z)}{\partial z_1}\left[\left(\frac{\partial F_i(z)}{\partial z_1}\right)^2+F_i(z)\frac{\partial^2 F_i(z)}{\partial z_1^2}\right]=4\frac{\partial \left(F_i(z)\frac{\partial F_i(z)}{\partial z_1}\right)^2}{\partial z_1}\eea
for $i=1,2$. Differentiating (\ref{md7}) partially with respect to $z_1$, we get
\bea\label{md9} \frac{\partial F_i(z+2c)}{\partial z_1}=-\left(\frac{\partial F_i(z)}{\partial z_1} \alpha_i(z)+\frac{\partial \alpha_i(z)}{\partial z_1} F_i(z)\right)\eea
for $i=1,2$. On the other hand, from (\ref{md2}) and (\ref{md8}), we obtain
\beas F_i(z)=\frac{2}{\alpha_i(z)}F_j(z+c)\frac{\partial F_j(z+c)}{\partial z_1}\eeas
and so (\ref{md2}) yields
\bea\label{md10} F_i(z)\frac{\partial F_i(z)}{\partial z_1}+\frac{1}{\alpha_j(z+c)}F_i(z+2c)\frac{\partial F_i(z+2c)}{\partial z_1}=0\eea
for $i,j\in\{1,2\}$ such that $i\neq j$. Consequently from (\ref{md7}), (\ref{md9}) and (\ref{md10}), we get
\beas F_i(z)\frac{\partial F_i(z)}{\partial z_1}+\frac{\alpha_i(z)}{\alpha_j(z+c)}F_i(z)\left(\frac{\partial F_i(z)}{\partial z_1} \alpha_i(z)+\frac{\partial \alpha_i(z)}{\partial z_1} F_i(z)\right)=0,\eeas
i.e.,
\bea\label{md11} \left(1+\frac{\alpha^2_i(z)}{\alpha_j(z+c)}\right)\frac{\partial F_i(z)}{\partial z_1}+\frac{\alpha_i(z)}{\alpha_j(z+c)}\frac{\partial \alpha_i(z)}{\partial z_1} F_i(z)=0\eea
for $i,j\in\{1,2\}$ such that $i\neq j$.
We divide following sub-cases.\par

\medskip
{\bf Sub-case 3.2.2.1.} Let $\frac{\partial \alpha_i(z)}{\partial z_1}\not\equiv 0$ for $i=1,2$. Note that $\parallel T(r,\alpha_i)=o(T(r,F_i))$ for $i=1,2$. Therefore from (\ref{md8}), it is easy to verify that
\bea\label{md12} \parallel\;N\left(r,0;\frac{\partial F_i(z)}{\partial z_1}\right)\leq N(r,0;\alpha_i)\leq T(r,\alpha_i)=o(T(r,F_i))\eea
for $i=1,2$. Consequently from (\ref{md11}) and (\ref{md12}), we conclude that
\bea\label{md13}\parallel\; N(r,0;F_i)=o(T(r,F_i))\eea
for $i=1,2$. Now (\ref{md8}) gives
\beas\frac{\alpha_i(z)}{F_i^3(z)}=8\left(\frac{\frac{\partial F_i(z)}{\partial z_1}}{F_i(z)}\right)^3+8\frac{\frac{\partial F_i(z)}{\partial z_1}}{F_i(z)}\times\frac{\frac{\partial^2 F_i(z)}{\partial z_1^2}}{F_i(z)}\eeas
and so by Lemma \ref{L.1}, we get
\bea\label{md14}\parallel\;m\left(r,0;F_i(z)\right)=o(T(r,F_i))\eea
for $i=1,2$. Now in view of the first main theorem and using (\ref{md13}) and (\ref{md14}), we obtain
$\parallel T(r,F_i)=o(T(r,F_i))$ for $i=1,2$, which is impossible.\par

\medskip
{\bf Sub-case 3.2.2.2.} Let $\frac{\partial \alpha_i(z)}{\partial z_1}\equiv 0$ for $i=1,2$. Now from (\ref{md11}), we get $\alpha^2_i(z)=-\alpha_j(z+c)$ for $i,j\in\{1,2\}$ such that $i\neq j$ and so
\bea\label{md15} \alpha^4_i(z)=-\alpha_i(z+2c)\eea
for $i\in\{1,2\}$.
If $\alpha_i(z)$ is transcendental, then using Lemma \ref{L.4} to (\ref{md15}), we get $\parallel m(r,\alpha_i)=o(T(r, \alpha_i))$ for $i=1,2$. Since $\alpha_i(z)$ are entire for $i=1,2$, we have $\parallel T(r, \alpha_i)=o(T(r, \alpha_i))$, leading to a contradiction. Hence $\alpha_i(z)$ is a polynomial in $\mathbb{C}^m$ for $i=1,2$.
It is easy to deduce from (\ref{md15}) that $\alpha_i(z)$ must be a constant, say $k_i(\neq 0)$ for $i=1,2$. Then from (\ref{md15}), we have $k_i^4=-k_i$, which implies that $k_i^3=-1$ for $i=1,2$. Since $\alpha^2_i(z)=-\alpha_j(z+c)$ for $i,j\in\{1,2\}$ such that $i\neq j$, we get $k_i^2=-k_j$ for $i,j\in\{1,2\}$ such that $i\neq j$.
Now differentiating (\ref{md8a}) partially with respect to $z_1$, we get
\bea\label{md16} 4\left(\frac{\partial F_i(z)}{\partial z_1}\right)^2\frac{\partial^2 F_i(z)}{\partial z_1^2}+F_i(z)\left( \left(\frac{\partial^2 F_i(z)}{\partial z_1^2}\right)^2+\frac{\partial F_i(z)}{\partial z_1}\frac{\partial^3 F_i(z)}{\partial z_1^3}\right)=0\eea
for $i=1,2$.
We consider following two sub-cases.\par

\medskip
{\bf Sub-case 3.2.2.2.1.} Let $\frac{\partial^2 F_i(z)}{\partial z_1^2}\not\equiv 0$ for $i=1,2$. Now in view of the first main theorem and using Lemma \ref{L.1}, we get
\bea\label{md17} \parallel\;N\left(r,0;\frac{\partial^2 F_i(z)}{\partial z_1^2}\right)&\leq& N\left(r,0;\frac{\frac{\partial^2 F_i(z)}{\partial z_1^2}}{\frac{\partial F_i(z)}{\partial z_1}}\right)+N\left(r,0;\frac{\partial F_i(z)}{\partial z_1}\right)\nonumber\\&\leq&
T\left(r,\frac{\frac{\partial^2 F_i(z)}{\partial z_1^2}}{\frac{\partial F_i(z)}{\partial z_1}}\right)+N\left(r,0;\frac{\partial F_i(z)}{\partial z_1}\right)+O(1)\nonumber\\&\leq&
N\left(r,\frac{\partial F_i(z)}{\partial z_1}\right)+2N\left(r,0;\frac{\partial F_i(z)}{\partial z_1}\right)+o(T(r,F_i))
\eea
for $i=1,2$. Since $F_i(z)$ is an entire functions, from (\ref{md12}) and (\ref{md17}), we obtain
\bea\label{md18}\parallel\;N\left(r,0;\frac{\partial^2 F_i(z)}{\partial z_1^2}\right)=o(T(r,F_i))\eea
for $i=1,2$.
Therefore using (\ref{md12}) and (\ref{md18}) to (\ref{md16}), we deduce that $\parallel N(r,0;F_i)=o(T(r,F_i))$ for $i=1,2$. Now proceeding in the same way as done in the proof of Sub-case 3.2.2.1, we get a contradiction.\par

\medskip
{\bf Sub-case 3.2.2.2.2.} Let $\frac{\partial^2 F_i(z)}{\partial z_1^2}\equiv 0$ for $i=1,2$. We know that $\alpha_i(z)=k_i$, where $k_i^3=-1$ and $k_i^2=-k_j$ for $i,j\in\{1,2\}$ such that $i\neq j$. Then from (\ref{md8a}), we obtain
\bea\label{md19} \frac{\partial F_i(z)}{\partial z_1}=\frac{K_i}{2},\eea
where $K_i^3=k_i$ for $i=1,2$. Clearly (\ref{md19}) implies
\bea\label{md20} \frac{\partial f_i(z)}{\partial z_1}=F_i(z)=\frac{K_i}{2}z_1+g_i(z_2,\ldots,z_m),\eea
where $g_i(z_2,z_3,\ldots, z_m)$ is a finite order transcendental entire function for $i=1,2$. On the other hand, from
(\ref{md7}), we have $F_i(z+2c)=-k_iF_i(z)$. Therefore (\ref{md20}) gives
\[\frac{K_i}{2}(z_1+2c_1)+g_i(z_2+2c_2, \ldots, z_m+2c_m)=-k_i\left(\frac{K_i}{2}z_1+g_i(z_2, z_3, \ldots, z_m)\right),\]
which implies $k_i=-1$ and so
\bea\label{md.21} g_i(z_2+2c_2, \ldots, z_m+2c_m)-g_i(z_2, \ldots, z_m)=K_ic_1\eea
for $i=1,2$. Again from (\ref{md20}), we have
\bea\label{md21} f_i(z)=\frac{K_i}{4}z_1^2 + z_1g_i(z_2,z_3,\ldots,z_m)+G_i(z_2,z_3,\ldots,z_m),\eea
where $G_i(z_2,z_3,\ldots, z_m)$ is a finite order transcendental entire function for $i=1,2$.
Since $F_i(z+2c)=F_i(z)$, from (\ref{md1}), we get $f_j(z+3c)=1-{F_i(z+2c)}^2=1-F_i^2(z)=f_j(z+c)$ for $i,j\in\{1,2\}$ such that $i\neq j$, which shows that $f_i(z)$ is $2c$-periodic for $i=1,2$. Then from (\ref {md21}), we get
\bea\label{md22} f_i(z)=f_i(z+2c)=\frac{K_i}{4}z_1^2 + z_1g_i(z_2,z_3,\ldots,z_m)+G_i(z_2,z_3,\ldots,z_m)\eea
for $i=1,2$.
Using (\ref{md20}) and (\ref{md22}) to (\ref{md3}), we get
\beas &&K_i^2\left(\frac{K_i^2}{4}z_1^2+K_iz_1g_i(z_2,\ldots,z_m)+g_i^2(z_2,\ldots,z_m)\right)\\&&+
\frac{K_i}{4}z_1^2 + z_1g_i(z_2,z_3,\ldots,z_m)+G_i(z_2,z_3,\ldots,z_m)=1,\eeas
i.e., $G_i(z_2,z_3,\ldots,z_m)=1-K_i^2(g_i(z_2,z_3,\ldots,z_m))^2$ for $i=1,2$
and so from (\ref{md22}), we have
\beas f_i(z)=1+\frac{K_i}{4}z_1^2 + z_1g_i(z_2,z_3,\ldots,z_m)-K_i^2g_i^2(z_2,z_3,\ldots,z_m),\eeas
where $K_i^3=-1$ and $g_i(z_2,z_3,\ldots, z_m)$ is a finite order transcendental entire function such that (\ref{md.21}) holds.\par

\medskip
{\bf Sub-case 3.3.} Let $m_2=n_2$. Then (\ref{lm.5}) gives $n_2=m_2=2$ and so from (\ref{pds1}), we get
\bea\label{md.1} F_2^2(z)+f_1^2(z+c)=1.\eea

Now using Theorem 1.2.A to (\ref{md.1}), we arrive at (\ref{llm.8})-(\ref{llm.10}).
Differentiating (\ref{lllm.1}) partially with respect to $z_1$, we get
\bea\label{md.5} n_1F_1^{n_1-1}(z+c)\frac{\partial F_1(z+c)}{\partial z_1}+F_2(z+2c)=0.\eea

Now using  (\ref{llm.8}) to (\ref{md.5}), we get
\bea\label{md.6}F_2^{n_1-1}(z)n_1\left(\frac{\partial h(z)}{\partial z_1}\right)^{n_1-1}\left(\frac{\partial^2 h(z)}{\partial z_1^2}F_2(z)+\frac{\partial h(z)}{\partial z_1}\frac{\partial F_2(z)}{\partial z_1}\right)=-F_2(z+2c).\eea

We know that the set of multiple zeros of $F_2(z)$ is algebraic. Therefore from (\ref{llm.8}) and (\ref{md.6}), we can easily conclude that $n_1=2$ and so
\beas 2\frac{\partial h(z)}{\partial z_1}\left(\frac{\partial^2 h(z)}{\partial z_1^2}\cos^2(h(z))-\frac{1}{2}\left(\frac{\partial h(z)}{\partial z_1}\right)^2\sin (2h(z))\right)=-\cos h((z+2c)).\eeas

Now proceeding in the same way as done in the proof of Sub-case 1.2, we get a contradiction.\par

Hence the proof is complete.
\end{proof}

\vspace{0.1in}
{\bf Compliance of Ethical Standards:}\par

{\bf Conflict of Interest.} The authors declare that there is no conflict of interest regarding the publication of this paper.\par

{\bf Data availability statement.} Data sharing not applicable to this article as no data sets were generated or analysed during the current study.

\end{document}